\documentclass[]{amsart}
\usepackage{amsmath, amsthm}

\usepackage{amssymb}
\usepackage{times}

\usepackage[all]{xypic}  \entrymodifiers={+!!<0pt,\fontdimen22\textfont2>}

\newtheorem{thm}{Theorem}[section]
\newtheorem{prop}[thm]{Proposition }
\newtheorem{defn}[thm]{Definition}
\newtheorem{lemma}[thm]{Lemma}
\newtheorem{cor}[thm]{Corollary}
\newtheorem{conj}[thm]{Conjecture}

\newtheorem{cond}[thm]{Condition}

\newcommand{\BM}{\operatorname{BM}}

\newcommand{\Ext}{\operatorname{Ext}}

\newcommand{\Spec}{\operatorname{Spec}}
\newcommand{\trdeg}{\operatorname{tr.deg}}
\newcommand{\im}{\operatorname{Im}}
\renewcommand{\ker}{\operatorname{ker}}
\newcommand{\coker}{\operatorname{coker}}

\newcommand{\colim}{\operatornamewithlimits{colim}}

\newcommand{\RHom}{\operatorname{RHom}}
\newcommand{\height}{\operatorname{ht}}

\newcommand{\calA}{\mathcal A}
\newcommand{\calB}{\mathcal B}
\newcommand{\calC}{\mathcal C}

\newcommand{\calF}{\mathcal F}
\newcommand{\calG}{\mathcal G}
\newcommand{\calH}{\mathcal H}

\newcommand{\calK}{\mathcal K}

\newcommand{\calO}{\mathcal O}

\newcommand{\calV}{\mathcal V}

\newcommand{\G}{{\mathbb G}}

\newcommand{\Q}{{{\mathbb Q}}}

\newcommand{\Z}{{{\mathbb{Z}}}}

\newcommand{\frakg}{\mathfrak g}

\newcommand{\frakp}{{\mathfrak p}}

\newcommand{\frakI}{{\mathfrak I}}

\numberwithin{equation}{subsection}

\title{Comparison of Dualizing Complexes}

\author{Changlong Zhong}

\begin{document}
\maketitle

\begin{abstract}
We prove that there is a map from Bloch's cycle complex to Kato's
complex of Milnor K-theory, which induces a quasi-isomorphism from
cycle complex mod $p^r$ to  Moser's complex of
logarithmic de Rham--Witt sheaves. Next we show that the truncation
of Bloch's cycle complex at $-3$ is quasi-isomorphic to Spiess'
dualizing complex.  In the end, we prove that a weak form of the Gersten Conjecture implies that Sato's dualizing complex is quasi-isomorphic to
Bloch's complex.

\end{abstract}

\section{Introduction} Using Lichtenbaum's weight-two
 motivic complex $\Z(X,2)$ (a two-term complex derived from relative K-theory),  M. Spiess \cite{spiess} constructed a complex of \text{\'{e}}tale sheaves $\calK_X$ (Definition \ref{definition:spiesscomplex}) on arithmetic surfaces $X$ over a Dedekind domain $D$ and used it to prove a duality theorem
   of constructible sheaves. For $X$ over a perfect field $k$ of characteristic $p>0$, T. Moser \cite{moser}  studied Gersten complexes of logarithmic de Rham--Witt sheaves $\widetilde{\nu}_{X,r}(n)$ (\ref{definition:Mosercomplex}) and showed that, when $k$ is finite, $\widetilde{\nu}_{X,r}\stackrel{def}=\widetilde{\nu}_{X,r}(0)$ is a dualizing complex for constructible $\Z/p^r$-sheaves. For regular semi-stable schemes $X$ over certain Dedekind domain $D$ (see Condition \ref{condition*}), K. Sato \cite{sato} defined certain dualizing complex $\frakI_r(n)_X$
($0\le n \le \dim X$) (Definition \ref{definition:satocomplex}) in derived category of \'{e}tale sheaves and proved a duality theorem for $\Z/p^r$-sheaves as well.
In more general
situations, for instance, schemes over algebraically closed fields,
finite fields, local fields and certain Dedekind domains, T. Geisser
\cite{geisser} proved that the complex $\Z^c_X\stackrel{def}=\Z^c_X(0)$
of \'{e}tale sheaves (see (\ref{definition:blochcomplex})) is a dualizing complex for constructible sheaves. Here $\Z^c_X(n)$ is Bloch's cycle complex whose homology defines higher Chow groups \cite{bloch}. In this paper, we answer the following questions of quasi-isomorphisms of these complexes:

\begin{thm}[Theorem \ref{thm:cycletoLRW}]
  For $X$ a scheme separated and essentially of finite type over a perfect field $k$ of characteristic $p>0$, and $n\le d=\dim X$, there is a map $$\hat{\psi}:\Z^c_X/p^r(n)\to \widetilde{\nu}_{X,r}(n)$$ which induces a quasi-isomorphism. Here $\Z^c_X/p^r(n)\stackrel{def}= \Z^c_X(n)/p^r$.
\end{thm}

\begin{thm}[Theorem \ref{thmcyclespiess} and Corollary \ref{corollaryRhom}]
 Let $X$ be a surface over a perfect field $k$ or a Dedekind domain $D$ with perfect residue fields. Then $\tau_{\ge -3}\Z^c_X$ is quasi-isomorphic to $\calK_X$. Moreover, for any torsion sheaf $\calF$, $$\RHom_X(\calF,\Z^c_X)\cong \RHom_X(\calF,\calK_X).$$
\end{thm}

\begin{thm}[Theorem \ref{thm:cycletosato2}]
 Let $p$ be a prime number, $X$ be a scheme over a Dedekind domain $D$ which satisfies Condition (\ref{condition*}) below. Assume that $d=\dim X$, and the conjecture $\calB(n)$ with $\Z/p^r$-coefficients (Conjecture \ref{conjecture:gersten}) holds for all $x\in X$. Then
$$\frakI_r(n)_X \stackrel{\cong}\longrightarrow (\Z^c_X/p^r(d-n)[-2d]).$$
\end{thm}

Let $X$ be a scheme of dimension $d$ over a perfect field of characteristic $p>0$. Moser's complex $\widetilde{\nu}_{X,r}(n)$ is the Gersten complex of logarithmic de Rham--Witt sheaf $W_r\Omega^{d-n}_{X,\log}$, when $X$ is smooth. For general $X$, $\widetilde{\nu}_{X,r}(n)(X)$ is identified with Kato's complex $C_X^M(n)/p^r$ of Milnor K-groups modulo $p^r$ (see Theorem \ref{thm:MKtoLRW}). To prove Theorem \ref{thm:cycletoLRW}, first, we show that the niveau filtration of higher Chow groups induces a canonical
map $$\phi:\Z^c_X(n)(X) \to C_X^{\text{HC}}(n),$$ the latter being the global sections of
Gersten complex of higher Chow groups (\'{e}tale sheafified). Then we show that
$C_X^{\text{HC}}(n)$ is isomorphic to  $C_X^{\text{M}}(n)$, which provides us a map
$$\psi:\Z^c_X(n)(X)\to C_X^{\text{M}}(n).$$ When composing with
the isomorphism $$C_X^M(n)/p^r\cong \widetilde{\nu}_{X,r}(n)(X),$$ we obtain a map
$$\hat{\psi}:\Z^c_X/p^r(n)=\Z^c_X(n)/p^r\to \widetilde{\nu}_{X,r}(n).$$ Moreover, using a result of
Geisser--Levine \cite[Theorem 1.1]{geisserlevine}, we show that $\hat{\psi}$ induces an
isomorphism of cohomology groups. Hence we
conclude that $\hat{\psi}$ is a quasi-isomorphism.  As another application of this method, we show that, for smooth and projective varieties over finite fields, the conjecture A($n$) of Geisser (part of  Parshin's Conjecture, see \cite[Proposition 2.1]{geisser2}) is true, if and only if $$\phi_{\Q}:\Q_X^c(n)(X) \to C_{X}^{HC}(n)\otimes \Q$$ is a quasi-isomorphism.

Let $X$ be a two dimensional scheme over a perfect field or a Dedekind domain with perfect residue fields. Spiess' complex $\calK_X$ is defined by connecting the complex $$\displaystyle\bigoplus_{X_{(1)}}i_{x*}\G_m \to \bigoplus_{X_{(0)}}i_{x*} \Z$$
with Lichtenbaum's weight-two motivic complex $\Z(k(X), 2)$ of the function field $k(X)$ (see \cite[Definition 2.1]{lichtenbaum} for $\Z(\_, 2)$). To compare it with Bloch's complex, we first define an intermediate complex $\calC(X)$. Using the niveau spectral sequence, we show that $\calC(X)$ is quasi-isomorphic to $\tau_{\ge -3}\Z^c_X(X)$. Then
we show that $\calC(X)$ is quasi-isomorphic to $\calK_X(X)$. A key
ingredient in this step is the quasi-isomorphism in \cite[\S7]{blochlichtenbaum} between a
truncation of Bloch's complex of $k(X)$ and $\Z(k(X),2)$. Based on the above, we conclude that $\tau_{\ge -3}\Z^c_X$
is quasi-isomorphic to $\calK_X$.

Let $X$ be a scheme over a Dedekind domain satisfying Condition (\ref{condition*}), and $\dim X=d$. Sato's complex $\frakI_r(n)_X$ is defined to be a complex in the derived category of \'{e}tale sheaves satisfying a distinguished triangle, with the other two terms the logarithmic de Rham-Witt sheaf on points of characteristic $p$ and $\mu_{p^r}^{\otimes n}$ on the open complement. In order to compare $\frakI_r(n)_X$ with $\Z^c_X/p^r(d-n)[-2d]$, we
show that $\Z^c_X/p^r(d-n)[-2d]$ satisfies a similar triangle (the localization
sequence of Bloch's cycle complex), and that there is a map between the two triangles, which induces isomorphisms on the cohomologies up to degree $n$. The key point is to show that these isomorphisms  are
compatible. Using the Gersten complexes of $\mu^{\otimes d-n}_{p^r}$, cycle complex and logarithmic de Rham--Witt sheaves, we replace
maps of cohomology of schemes by those of function fields, whose compatibility is
straightforward. Hence the two complexes are quasi-isomorphic up to degree $n$. With the aid of the conjecture
$\calB(n)$ with $\Z/p^r$-coefficients, then they are actually quasi-isomorphic.

 The paper is organized as follows: in Section 2, we recall the definitions of Bloch's cycle complex $\Z^c_X(n)$, Kato's complex $C_X^M(n)$ of Milnor K-theory and
 Moser's complex $\widetilde{\nu}_{r}(n)$, as well as the duality results of Geisser and Moser. We also recall the construction of the niveau spectral sequence
 of higher Chow groups. For the reader's convenience, we include here a generalization of the Beilinson--Lichtenbaum Conjecture to general schemes and recall Levine's proof of Kummer isomorphism for regular schemes over Dedekind domains  \cite[Levine, Theorem 12.5]{levine2}, with the assumption of the conjecture $\calB(n)$. In Section 3, first we recall Spiess complex $\calK_X$, his duality results and the definition of Lichtenbaum's weight-two motivic complex $\Z(X,2)$. Then we define the complex $\calC(X)$ and  compare it with $\Z^c_X$ and $\calK_X$, respectively. In Section 4, we recall Sato's definition of $\frakI_r(n)_X$ ($0\le n \le \dim X$), and prove Proposition \ref{prop:cycletosato} and Theorem \ref{thm:cycletosato2}.

\section*{Terminology} Throughout this paper, the concepts
chain complex and cochain complex are used interchangeably. For
instance, if $A$ is a chain complex, we think of it as a cochain
complex by letting $A^n=A_{-n}$. The convention for shift is:
$A[n]^i=A^{i+n}$. When talking about truncations, we mean truncations in cohomological degrees.

We use $D$ to denote a Dedekind domain of
characteristic $0$ with perfect residue fields, and $k$ to denote a
perfect field of characteristic $p\ge 0$. There is special assumption for $D$ in Section 4 (Condition \ref{condition*}). All the schemes in this paper will
be separated and essentially of finite type over $S$ with $S=\Spec
k$ or $\Spec D$. By variety we mean schemes separated and of finite type over fields.
The
dimension of an irreducible $S$-scheme $X$ (or dimension of $X$ over $S$) is
defined as
 $$\dim_SX\stackrel{def}=\trdeg(k(X):k(\frakp))-\height \frakp+\dim S$$
  where $\frakp$ is the image of the generic point of $X$ in $S$, and $k(X)$ is the function field of $X$. If $X$ is of finite type over $S$, then $\dim_SX=\dim X$, the Krull dimension of $X$. If $X$ is spectrum of a field of transcendental degree $m$ over $S=\Spec k$, then $\dim_SX=m$. Note that this definition of dimension is different from that of relative dimensions. For instance, if $S=\Spec D$ and $X$ has relative dimension $d$ over $S$, then $\dim_SX=d+1.$
  We use $X_{(i)}$ (resp. $X^{(i)}$) to denote the set of dimension $i$ (resp. codimension $i$) points of $X$. For $x\in X$, we use $i_x$ to denote the embedding of $x$ and $i_{x*}$ the push-forward of (\text{\'{e}}tale) sheaves.
  The number  $n$ is always less than $d=\dim X$.

\section*{Acknowledgements}
I am grateful to my advisor Thomas Geisser for his instruction and
help. I am indebted to the referee for pointing out a mistake in a previous version. I also want to thank Enxin Wu for helpful conversations.

\section{Comparison between Bloch's Complex and Moser's Complex}

\noindent{\textit{Bloch's cycle complex.}}~~ Let
$$\Delta^i=\Delta^i_S=S\times_\Z\Spec \Z[t_0,...,t_i]/(\sum t_j-1).$$
We define $z_n(X,i)$ as the free abelian group generated by closed integral subschemes
$Z\subset X\times \Delta^i$ that intersect all the faces properly
and $\dim_SZ=n+i$. Then $z_n(\_,i)$ is a sheaf
in the Zariski and \text{\'{e}}tale topology on $X$. The following complex is defined
by Bloch (\cite{bloch}, or see \cite{geisser} for the
notations):
\begin{equation}
\to z_n(\_,i) \stackrel{d}\rightarrow ...\to z_n(\_,1) \to z_n(\_,0)
\to 0,\label{definition:blochcomplex}
\end{equation}
where
$$d(Z)=\sum_j(-1)^j[Z\cap V(t_j)],$$ with $V(t_j)$ the closed integral subscheme generated by $t_j$ and $[Z\cap V(t_j)]$ the  linear combination
of irreducible components of $Z\cap V(t_j)$ with coefficients
intersection multiplicities. We define the complex $\Z^c_X(n)^t$ to be the complex
of sheaves in the topology $t$ with $t$=Zar or \'{e}t, and put $z_n(\_,-i-2n)$ in
(cohomological) degree $i$, i.e.,
$$(\Z^c_X(n)^t)^i=z_n(\_,-i-2n).$$
 In this paper we are mostly interested in
\'{e}tale sheaves, so by $\Z^c_X(n)$ we  mean the complex of \'{e}tale sheaves,
unless otherwise stated. Note that $\Z^c_X(n)^{Zar}$ satisfies Zariski descent. Since $\Z^c_X(n)$ is a
complex of flat sheaves, the derived tensor product agrees
with the usual tensor product of complexes. For an abelian group $A$, define
$$A_X^{c}(n)=A\otimes \Z_X^c(n).$$  The complex $\Z^c_X(n)(X)$ is covariant for proper maps and
contravariant for quasi-finite, flat maps.
Let $\Z^c_X=\Z^c_X(0)$, and omit $X$ if there is no ambiguity. Define motivic Borel--Moore homology
to be
$$H_i^{\BM}(X/S,\Z(n))\stackrel{def}=H_i(\Z^c_X(n)(X))=H^{-i}(\Z^c_X(n)(X)).$$
This definition of $\Z^c_X(n)$ and $H^{\BM}_i(X/S, \Z(n))$ depend on the base scheme $S$. Unless otherwise specified, they will be defined over base $S$ and we will just write $\Z^c_X(n)$ and $H_i^{\BM}(X,\Z(n))$, respectively. As an example, if $F$ has transcendental degree $d$ over the base $k$, then $$H^{\BM}_{2d+i}(F/k,\Z(d+n))= H_{i}^{\BM}(F/F,\Z(n)).$$
For a scheme $X$ of pure dimension $d$ over $S$,
$$H_i^{\BM}(X,\Z(n))=CH^{d-n}(X,i-2n).$$
For a scheme $X$ smooth over
 over a perfect field $k$ and $\dim X=d$,
 $$\Z^c_X(d-n)^{Zar}\cong\Z(n)[2d].$$
 Here $\Z(n)$ is the motivic complex defined by Voevodsky, and the ~$\cong$~ here is a quasi-isomorphism. In particular, for a field $F$ of
transcendental degree $d$ over $S=\Spec k$,
$$H_i^{\BM}(F/k,\Z(n))=CH^{d-n}(F,i-2n)\cong H^{2d-i}(F,\Z(d-n)).$$

The following three theorem state important property of cycle complex, and we will use them in Section 4. The first one is  the localization property of cycle complex.

\begin{thm}[Levine, \mbox{\cite[Theorem 0.7]{levine}}]\label{thm:localization}
  For any $X$, let $\xymatrix{U \ar@{^(->}[r]^j & X}$ be an open subscheme with closed complement $\xymatrix{Y\ar@{^(->}[r]^i & X}$. Then there is a distinguished triangle of complexes of Zariski sheaves
  $$\xymatrix{
  \ar[r] & i_*\Z^c_Y(n)^{Zar} \ar[r] & \Z^c_X(n)^{Zar} \ar[r] & j_*\Z^c_U(n)^{Zar} \ar[r] &.
  }$$
\end{thm}

We will need a weak form of the Gersten Conjecture of higher Chow groups.
\begin{conj}[Conjecture $\mathcal{B}(n)$]\label{conjecture:gersten}
Let $R$ be a regular local ring and $K$ be its function field.  Then the canonical map
$$H^{\BM}_{s}(R,\Z/m(n))\to H^{\BM}_s(K,\Z/m(n))$$
is injective for any $s$.
\end{conj}

\begin{lemma}\label{lemma:gerstenconjecture}
If $X$ is regular with dimension $d$, and the conjecture $\calB(n)$ is true for all the local rings of points $x\in X$, then the complexes $\Z^c_X/m(n)$ and $\Z^c_X/m(n)^{Zar}$ are acyclic at degree $>-d-n$.
\end{lemma}

\begin{proof}
  From the conjecture, we see that for any $x\in X$ with local ring $\calO_{X,x}$ and function field $K_x$, the maps are injective
  $$H^{\BM}_s(\calO_{X,x},\Z/m(n))\to H^{\BM}_s(K_x,\Z/m(n)).$$
  But $$H^{\BM}_{s}(K_x, \Z/m(n))=H^{2d-s}(K_x, \Z/m(d-n)).$$
   Here $H^{2d-s}(K_x,\Z/m(d-n)$ is Voevodsky's motivic cohomology, which vanishes for $s<d+n$ \cite[Theorem 3.6]{voevodsky}. Hence,  $H^{\BM}_s(\calO_{X,x}, \Z^c/m(n))=0$ for $s<d+n$. It implies that $\calH^i(\Z^c_X/m(n)^{Zar})=0$ for $i>-d-n$. The same property holds for $\Z^c_X/m(n)$. Hence we finish the proof.
\end{proof}

Let $\epsilon$ be the map from \'{e}tale site to Zariski site. The Beilinson--Lichtenbaum Conjecture claims that the following map
$$\eta_X:\Z^c_X/m(n)^{Zar}\to\tau_{\le -d-n}R\epsilon_*\Z^c_X/m(n)$$
is an quasi-isomorphism. Recall that the Bloch--Kato Conjecture claims that for a field $F$ and $m\in \Z$, the Galois symbol
$$K_n^M(F)/m\to H^n(F_{\acute{e}t}, \mu^{\otimes n}_m)$$
is an isomorphism.  Since it was proved by Voevodsky and Rost recently, we will refer it as the Rost--Voevodsky Theorem. It implies the Beilinson--Lichtenbaum Conjecture for smooth varieties over a field (see \cite{suslinvoevodsky} if the field has characteristic 0, and \cite{geisserlevine2} for positive characteristic case). We generalize it to general schemes over $S$, assuming the conjecture $\calB(n)$. First we prove a lemma from homological algebra.

\begin{lemma}
  Let $F:\calA \to \calB$ be a left exact functor, and assume that there are enough injectives in $\calA$. Assume that $C^*$ is a (cochain) complex. Then $$\tau_{\le n}RF(\tau_{\le n}C^*)\cong \tau_{\le n}RFC^*.$$
\end{lemma}
\begin{proof}
Let $C^*\to I^*$ be an injective resolution, i.e., $I^*$ is a complex of injectives and is quasi-isomorphic to $C^*$. Let $A$ be the kernel of $$d_n:I^n\to I^{n+1}.$$ Then we obtain an injective resolution of $A$ as follows:
$$\xymatrix{
0 \ar[r]& A \ar@{^(->}[r] &I^n\ar[r]^{d_n} &I^{n+1} \ar[r]^{c_{n+1}} & J^{n+1} \ar[r] & .
}$$
Then the new complex
$$\xymatrix{
\ar[r] & I^{n-1} \ar[r]^{d_{n-1}} & I^n\ar[r]^{d_n} & I^{n+1} \ar[r]^{c_{n+1}} & J^{n+1} \ar[r] & } $$
is an injective resolution of $\tau_{\le n} C^*$, which we denote as $J^*$. Then for $m\le n$, $$H^m(FI^*)\cong H^m(FJ^*).$$ Hence we get the conclusion.

\end{proof}

\begin{thm} Let $X$ be a pure dimensional scheme over $S$ with $S=\Spec k$ or $\Spec D$, and $\dim_SX=d$.   Then there is an quasi-isomorphism
$$\eta_X:\tau_{\le -d-n}\Z^c_X/m(n)^{Zar}\to\tau_{\le -d-n}R\epsilon_*\Z^c_X/m(n).$$
\label{thm:BLconjecturetorsion}
If $X$ is regular and the conjecture $\calB(n)$ with $\Z/m$-coefficients is true for all the points $x\in X$, then the truncation in front of the first item can be removed.

\end{thm}
\begin{proof} If $X$ is smooth over $S=\Spec k$, then the conclusion $$\Z^c_X/m(n)^{Zar}\to\tau_{\le -d-n}R\epsilon_*\Z^c_X/m(n)$$  is implied by the Rost-Voevodsky Theorem. If $X$ is a general variety over $\Spec k$ we prove the conjecture by induction on the dimension
of $X$. Suppose that it is true for $Z$ such that $\dim Z<d$. Now suppose $\dim X=d$. We can assume that $X$ is reduced. Let $m=-d-n$ and $U\stackrel{j}\longrightarrow X$ be an smooth open subscheme with $Z\stackrel{i}\longrightarrow X$ its complement. Enlarge $Z$ so that it has pure codimension $1$. From Theorem \ref{thm:localization}, there is a distinguished triangle
  $$\xymatrix{
  \ar[r] & i_*\Z^c_Z/m(n)^{Zar} \ar[r] & \Z^c_X/m(n)^{Zar} \ar[r] & j_*\Z^c_U/m(n)^{Zar} \ar[r] &
  .}$$
  Apply the functor $\epsilon^*$, then apply the functor $R\epsilon_*$, we obtain a distinguished triangle
$$\xymatrix{R\epsilon_*(\epsilon^*i_*\Z^c_Z/m(n)^{Zar}) \ar[r] & R\epsilon_*(\epsilon^*\Z^c_X/m(n)^{Zar})\ar[r] & R\epsilon_*(\epsilon^*j_*\Z^c_U/m(n)^{Zar}) .}$$
Note that \cite[Proposition 2.2]{geisser3}
$$\epsilon^*\Z^c_X/m(n)^{Zar}=\Z^c_X/m(n), \qquad R\epsilon_*(\epsilon^*i_*\Z^c_Z/m(n)^{Zar})=i_*R\epsilon_*(\Z^c_Z/m(n)).$$ By induction assumption, $$\tau_{\le -(d-1)-n}i_*R\epsilon_*(\Z^c_Z/m(n))\cong \tau_{\le -(d-1)-n}i_*\Z^c_Z/m(n)^{Zar}.$$ Hence, by the five lemma, it suffices to show that $$\tau_{\le m}R\epsilon_*(\epsilon^*j_*\Z^c_U/m(n)^{Zar}) \cong \tau_{\le m}j_*\Z^c_U/m(n)^{Zar}.$$ That is, we only have to compare
$R^i\epsilon_*(\epsilon^*j_*\Z^c_U/m(n)^{Zar})$ with $H^i(j_*\Z^c_U/m(n)^{Zar})$ for $i\le m$. The former one is the Zariski sheaf associated to the presheaf sending $V$ to $$H^i_{et}(V,\epsilon^*j_*\Z^c_U/m(n)^{Zar}),$$ while the latter one is the Zariski sheaf associated to the presheaf sending $V$ to (note that $\Z^c/m(n)^{Zar}$ satisfies Zariski descent) $$H^i_{Zar}(V\times_XU, \Z^c_U/m(n)^{Zar}).$$  Since $U$ is smooth, $\Z^c_U/m(n)^{Zar}$ satisfies the Beilinson--Lichtenbaum Conjecture, so the two sheaves are isomorphic. Hence,   $$\tau_{\le -d-n}\Z^c_X/m(n)^{Zar}\to\tau_{\le -d-n}R\epsilon_*\Z^c_X/m(n).$$

If $X$ is a smooth scheme over  $S=\Spec D$ with perfect residue fields, then the isomorphism $$\Z^c_X/m(n)^{Zar}\to\tau_{\le -d-n}R\epsilon_*\Z^c_X/m(n)$$ is proved in \cite[Theorem 1.2(2)]{geisser3}. Now let $X$ be a scheme flat over $D$. Since this is a local problem, we can assume that $D$ is a discrete valuation ring. Let $U$ be the generic fiber and  $Z$ be the special fiber. Then $\dim_SZ=d-1$. So we have a localization diagram similar as above. Both $U$ and $Z$ are varieties over fields, and $U$ is smooth, so similarly, we can prove the first part.

If the conjecture $\calB(n)$ with $\Z/m$-coefficients is true, from Lemma \ref{lemma:gerstenconjecture}, we get the conclusion.

 \end{proof}

For the reader's convenience, we recall Levine's proof of a weak form of the Kummer isomorphism for regular schemes over Dedekind domains.

\begin{thm}[Levine, \mbox{\cite[Theorem 12.5]{levine2}}]\label{thm:kummerisomorphism}
If $X$ is a regular and  pure-dimension $d$ scheme over $S=\Spec D$ or $\Spec k$, $m$ is invertible in $S$. Then the following map is an isomorphism
$$\tau_{\le -d-n}\Z^c_X/m(n)\cong \mu_{m}^{\otimes d-n}[2d].$$
If the conjecture $\calB(n)$ with $\Z/m$-coefficients holds for all the points $x\in X$, then the truncation can be removed.
 \end{thm}
 \begin{proof} If $X$ is smooth over a field, then it is proved in \cite[Theorem 1.5]{geisserlevine2}. If $X$ is smooth over $D$, then it is proved in \cite[Theorem 1.2(4)]{geisser3} (in this case the map $\varrho_X:\Z^c_X/m(n)\to \mu_m^{\otimes d-n}[2d]$ is defined in \cite[\S12]{levine2}).   The only case left is when $X$ is a regular scheme flat over $D$, so $d\ge 1$. This is a local problem, so assume that $D$ is a discrete valuation ring with perfect residue field of characteristic coprime to $m$.

    First we claim that, if $X$ is regular and admits an closed-open decomposition
$$X=Z\coprod U$$
such that $Z$ smooth and the conclusion holds for $U$, then it is also true for $X$.

By localization sequence of cycle complex and $\mu_m$, there is a commutative diagram
$$\xymatrix{
i_*\Z^c_Z/m(n) \ar[d] \ar[r] & \Z^c_X/m(n)\ar[d] \ar[r]& \epsilon^*j_*\Z^c_U/m(n)^{Zar}\ar[d]\\
i_*\mu^{\otimes b-n}_{m,Z}[2b] \ar[r] & \mu^{\otimes d-n}_{m,X}[2d] \ar[r]& Rj_*\mu^{\otimes d-n}_{m,U}[2d]
.}$$
Here the two rows are distinguished triangles, and $b=\dim Z$. The existence of the lower distinguished triangle is due to Gabber, \cite{fujiwara}. Then
$$R^{s}(j_*\mu_{m,U}^{\otimes d-n}[2d])=\left\{\begin{array}{cc} \mu_{m,X}^{\otimes d-n}[2d], & s=-2d,\\
                                  i_*\mu^{b-n}_{m,Z}[2b], & s=-1-2b,\\
                                  0,  &\textrm{otherwise}.
                                  \end{array}\right.$$
Moreover, we have
\begin{equation}\tau_{\le -d-n}\epsilon^*j_*\Z^c_U/m(n)^{Zar}\cong \tau_{\le -d-n}Rj_*\Z^c_U/m(n) \cong\tau_{\le -d-n}(Rj_*\mu^{\otimes d-n}_{m,U}[2d]).\end{equation}
Here the first isomorphism follows from the Beilinson--Lichtenbaum Conjecture for regular $U$ (Theorem \ref{thm:BLconjecturetorsion}), and the second one follows from the assumption on $U$.

If $n\le b$, then since $Z$  is smooth, we have isomorphism
$$\Z^c_Z/m(n)\cong \mu^{\otimes b-n}_{m,Z}[2b].$$ Hence, the two distinguished triangles have  isomorphic cohomology up to degree $-d-n$. In particular,
$$\tau_{\le -d-n}\Z^c_X/m(n)\cong \tau_{\le -d-n}(\mu_{m,X}^{\otimes d-n}[2d])\cong \mu_{m,X}^{\otimes d-n}[2d].$$

If $n>b$, then $\Z^c_Z/m(n)=0$. Then
$$\tau_{\le -d-n}\Z^c_X/m(n)\cong \tau_{\le -d-n}\epsilon^*j_*\Z^c_U/m(n)^{Zar}\cong \tau_{\le -d-n}(Rj_*\mu_{m,U}^{\otimes d-n}[2d]).$$
But $b\le d-1 \textrm{~and ~} n>b$ imply that $-2d\le -d-n<-1-2b.$ So $$\tau_{\le -d-n}(Rj_*\mu_{m,U}^{\otimes d-n}[2d])\cong R^{-2d}(j_*\mu_{m,U}^{\otimes d-n}[2d])\cong \mu_{m,X}^{\otimes d-n}[2d].$$
This finishes the proof of the claim.

Now for $X$ regular over $D$, we can find $U$ and a filtration of $X$
$$\emptyset =Z_{{t+1}}\subset Z_{{t}}\subset... \subset Z_{1}\subset Z_{0}=X $$
such that $Z_1$ and $U$ are the special and generic fiber of $X$, respectively, and $Z_{{i+1}}$ is the singular locus of $Z_{i}$. Then we have the following decompositions
$$X-Z_{2}=(Z_{1}-Z_{2})\coprod (X-Z_{1}),$$
$$X-Z_{3}=(Z_{2}-Z_{3})\coprod (X-Z_{2}),$$
$$...$$
$$X=X-Z_{{t+1}}=(Z_{{t}}-Z_{{t+1}})\coprod (X-Z_{t})$$
Here all the $Z_{i}-Z_{{i+1}}$ and $X-Z_1=U$ are smooth, and all the $X-Z_{i}$ are regular. From the first decomposition, we get the existence of the isomorphism for $X-Z_2$. Inductively, we get that the conclusion is true for $X=X-Z_{t+1}$.

If the conjecture $\calB(n)$ with $\Z/m$-coefficients is true, then the truncation on $\Z^c_X/m(n)$ is not necessary.
\end{proof}

The following are the duality results by Geisser, \cite{geisser}.

\begin{thm}[Geisser, \mbox{\cite[Theorem 4.1]{geisser}}] If $f:X \to \Spec k$ is separated and of finite type over a perfect field, then for any constructible sheaf $\calF$ on $X$, there is a canonical quasi-isomorphism
$$\RHom_X(\calF, \Z^c_X)\cong \RHom_{\Spec k}(Rf_!\calF, \Z).$$
\end{thm}
 \begin{cor}[Geisser, \mbox{\cite[\S5]{geisser}}](1) Let $S$ be the spectrum of a number ring or a finite field, and $\calF$ be a constructible sheaf on $X$, then there are perfect pairings of finite groups
$$H^i_c(X_{\text{\text{\'{e}}t}},\calF) \times \Ext^{2-i}_X(\calF, \Z^c_X)\to \Q/\Z.$$

(2)Let $k$ be an algebraically closed field and $\calF$ be a
constructible sheaf on $X$. Then there are perfect pairings of
finitely generated groups
$$H^i_{c}(X_{\text{\text{\'{e}}t}},\calF) \times \Ext^{1-i}_X(\calF, \Z^c_X) \to \Q/\Z.$$
\end{cor}

\noindent\textit{Niveau Spectral Sequence.}~~
The method of this section depends highly on the niveau
spectral sequence of higher Chow groups, so let us recall its construction.
We
adopt the notation in \cite[Proposition 2.1]{geisser}. We only define it for $\Z$-coefficients. Such spectral
sequence with other coefficients is totally analogous. Let $p$ be
the projection $X\times \Delta^i\to X$. Let  $F_s=F_s\Z^c(n)(X)$ be
the subcomplex generated by cycles $Z\subset X\times \Delta^i$ with
$\dim p(Z)\le s$, so $$(F_s)_m=(F_s)^{-m}\subset
(\Z^c(n)(X))^{-m}=z_n(X,m-2n).$$ There is a short exact sequence of
complexes:
$$
 0 \to F_{s-1}\to F_s \to F_s/F_{s-1} \to 0,
$$
which induces a long exact sequence of abelian groups
\begin{equation}
 \to H_{s+t+1}(F_{s}/F_{s-1}) \to H_{s+t}(F_{s-1}) \to
 H_{s+t}(F_s) \to H_{s+t}(F_s/F_{s-1}) \to \label{section 4: longexactsequence}.
\end{equation}
Moreover, by the localization property of higher Chow groups, we
have
$$
H_{s+t}(F_s/F_{s-1})\cong \bigoplus_{X_{(s)}} H_{s+t}^{\BM}(k(x),\Z(n))
.
$$
Therefore, there is a convergent spectral sequence:
\begin{equation}E^1_{s,t}=\displaystyle\bigoplus_{X_{(s)}} H_{s+t}^{\BM}(k(x),\Z(n))\Rightarrow H_{s+t}^{\BM}(X,\Z(n)).\label{spectralsequence:niveau}\end{equation}
This construction induces a filtration
$$N_sH_{m}^{\BM}(X,\Z(n))\stackrel{def}=\im (H_m(F_s)\to H_m^{\BM}(X,\Z(n))),$$
which is called the niveau
filtration of higher Chow groups. Note that $E^{1}_{s,t}=0$ for
$t<n$.

From the niveau spectral sequence (\ref{spectralsequence:niveau}), we can define $C_X^{\text{HC}}(n)\stackrel{def}=E^1_{*,n}$
 with $$(C_X^{\text{HC}}(n))^i=\bigoplus_{X_{(-i-n)}}H_{-i}^{\BM}(k(x), \Z(n))=\bigoplus_{X_{(-i-n)}}CH^{-i-2n}(k(x), -i-2n).$$
   and  $C^{\text{HC}}_X(n)/p^r$ to be the corresponding complex with terms $$\bigoplus_{X_{(-i-n)}}H_{-i}^{\BM}(k(x),\Z/p^r(n)).$$ If $X$ is a scheme over a perfect field of characteristic $p>0$, then $$H_{-i}^{\BM}(k(x), \Z/p^r(n))\cong H_{-i}^{\BM}(k(x), \Z(n))/p^r,$$ and $C^{\text{HC}}_X(n)/p^r$ is exactly the complex $C^{\text{HC}}_X(n)$ modulo $p^r$, so in this case our notation causes no ambiguity. Denote by
$(C_X^{\text{HC}}(n))^{t}$ the
sheafification in topology $t$ of the complex of presheaves $C_X^{\text{HC}}(n)$.
 This complex $(C_X^{\text{HC}}(n))^t$ is called the Gersten complex of  sheaf $H_{d+n}(\Z^c_X(n)^t)$. It is called  the Gersten
resolution of $H_{d+n}(\Z^c_X(n)^t)$ if they are quasi-isomorphic via
the canonical map $$H_{d+n}(\Z^c_X(n)^t)\to
\bigoplus_{X_{(d)}}a_ti_{x,p}H_{d+n}^{\BM}(\_, \Z(n)).$$
Here $a_t$ is denoted as the sheafification in topology $t$, and $i_{x,p}$ is the direct image of presheaves along the embedding $i_x:x\to X$. Note that $a_ti_{x,p}\neq i_{x*}a_t$ if $t=\acute{e}t$ and $i_{x*}$ is the direct image of \'{e}tale sheaves.

There is a map of complexes $$\phi:\Z^c_X(n)(X)\to
C_X^{\text{HC}}(n)$$ which induces the edge morphisms of the spectral
sequence. More precisely, $$\phi_m:(\Z^c_X(n)(X))_{m}=z_n(X,m-2n)
\to \bigoplus_{X_{(m-n)}}H^{\BM}_{m}(k(x),
\Z(n))=(C_X^{\text{HC}}(n))_{m}$$ satisfies the following
properties:

i) if $m>n+d$, then $m-n>d$ and $\phi_m=0$; similarly, if $m<n$,
$\phi_m=0$ as well.

ii) if $n\le m\le n+d$, then $\phi_m$ is the composition:
\begin{equation}z_n(X,m-2n)=(F_{m-n})_{m}\twoheadrightarrow \frac{(F_{m-n})_{m}}{(F_{m-n-1})_{m}} \overset{\cong}{\longrightarrow} \bigoplus_{X_{(m-n)}}z_n(k(x),m-2n) \twoheadrightarrow \bigoplus_{X_{(m-n)}}H_{m}^{\BM}(k(x),\Z(n)).
\label{filtration}\end{equation}

\noindent\textit{Gersten complex of logarithmic de Rham--Witt
sheaves.}~~
Let $X$ be a scheme essentially of finite type over $k$ of characteristic $p>0$. If $m\ge 0$, then
define the logarithmic de Rham--Witt sheaf $W_r\Omega^m_{X,\log}$ to be the \'{e}tale subsheaf of $W_r\Omega^m_{X}$ generated by elements of the form
$$d\log f_1\wedge ...\wedge d\log f_m$$
with $f_i\in \calO_X^*$. Define $W_r\Omega^m_{X,\log}=0$ if $m<0$. If $X=\Spec F$ is a field, we use the notation $\nu_{F,r}^m=W_r\Omega^m_{X,\log}$. If $X$ is smooth, define $\widetilde{\nu}_{X,r}(n)$ to be the Gersten
complex of logarithmic de Rham--Witt sheaves:
\begin{equation}
 0 \to \bigoplus_{X_{(d)}}i_{x*}\nu^{d-n}_{k(x),r} \to ... \to \bigoplus_{X_{(1)}} i_{x*}\nu^{1-n}_{k(x),r} \to \bigoplus_{X_{(0)}} i_{x*} \nu^{-n}_{k(x),r}  \to 0. \label{definition:Mosercomplex}
\end{equation}
The term indexed by $X_{(-i-n)}$ is in  (cohomological) degree $i$, i.e.,
 $$\widetilde{\nu}_{X,r}(n)^i=\bigoplus_{X_{(-i-n)}}\nu^{-i-2n}_{k(x),r}.$$
 Define $\widetilde{\nu}_{X,r}=\widetilde{\nu}_{X,r}(0).$ The
differentials of $\widetilde{\nu}_{X,r}(n)$ are induced from the coniveau
filtration of the sheaf $W_r\Omega^{d-n}_{X,\log}$, similar as the one for higher Chow groups above.
 If $X$ is not smooth, using
Theorem \ref{thm:MKtoLRW} below, Moser \cite{moser} also defines $\widetilde{\nu}_{X,r}(n)(X)$ by identifying it with $C_X^{\text{M}}(n)/p^r$ (see definition below). The following is the duality result of Moser.

\begin{thm}[Moser, \mbox{\cite[Theorem 5.6]{moser}}]
Let $k$ be a finite field of characteristic $p$ and $X$ be a
$k$-scheme of pure dimension $d$. Then for every $r\ge 1$ and every
constructible $\Z/p^r$-sheaf $\calF$, there are perfect pairings of
finite groups
$$H^i_c(X_{\text{\text{\'{e}}t}},\calF) \times \Ext_X^{1-i}(\calF, \widetilde{\nu}_{X,r}) \rightarrow \Z/p^r.$$
\end{thm}

\noindent\textit{Kato's complex of Milnor K-theory.}~~Let $X$ be a
scheme over $S$. Define $C_X^{\text{M}}(n)$ to be Kato's complex of
Milnor K-theory (cf. K. Kato, \cite{kato1}):
\begin{equation}\label{definition:katocomplex}
0 \to \bigoplus_{X_{(d)}} K_{d-n}^M(k(x)) \xrightarrow{d'} ...\to
\bigoplus_{X_{(1)}}K_{1-n}^M(k(x)) \to
\bigoplus_{X_{(0)}}K_{-n}^M(k(x)) \to 0.
\end{equation}
The group $K_m^M(k(x))$ is Milnor K-group, and $K_m^M(k(x))=0$ if
$m<0$. The differential $d'$ is defined as follows: for any $x\in
X_{(m)}$ and any $y\in \overline{\{x\}}\cap X_{(m-1)}$, we take the
normalization $\overline{\{x'\}}$ of $\overline{\{x\}}$ with $x'$
its generic point and define a map
$$\partial_y:K_{m-n}^M(k(x))=K_{m-n}^M(k(x'))\stackrel{\sum
\partial_{y'}}\longrightarrow \bigoplus_{y'|y}K^M_{m-n-1}(k(y'))
\stackrel{\sum N_{k(y')/k(y)}}\longrightarrow K_{m-n-1}^M(k(y)).$$
Here the notation $y'|y$ means that $y'\in
\overline{\{x'\}}_{(m-1)}$ is in the fiber of $y$,
$$N_{k(y')/k(y)}:K_{m-n-1}^M(k(y')) \to K_{m-n-1}^M(k(y))$$ is the
norm map of Milnor K-theory (see Bass--Tate, \cite{basstate} and
Kato, \cite[Section 1.7]{kato}), and $\partial_{y'}$ is the tame
symbol defined by $y'$. Then
$$d'\stackrel{def}=\displaystyle\sum_{y\in X_{(m-1)}\cap
\overline{\{x\}}} \partial_y.$$ Note that the sum in $d'$ is
finite since elements in $K_{m-n}^M(k(x'))$ are represented by sums of tensors of $m-n$
elements
 in $k(x')^*$, and each element in $ k(x')^*$ has a finite number of poles and zeros. When applying the tame symbol,
 only a finite number of terms in the sum are non-zero. Suppose that in (cohomological) degree $i$,
 $$C_X^{\text{M}}(n)^i=\bigoplus_{X_{(-i-n)}}K_{-i-2n}^M(k(x)),$$
  and  define $C_X^{\text{M}}=C_X^{\text{M}}(0)$. The complex $C_X^{\text{M}}(n)$ is covariant for proper maps and contravariant for quasi-finite and flat maps (see Rost, \cite[Proposition 4.6(1),(2)]{rost}).

  Note that all the three complexes, $C_X^{\text{HC}}(n), C_X^{\text{M}}(n) $ and $\widetilde{\nu}_{X,r}(n)$ are concentrated in cohomological degree
  $[-n-d, -2n]$ when  $ n\le d=\dim X$.

\begin{thm}[Bloch--Kato, Gros--Suwa, Moser]
For $X$ a scheme separated and essentially of finite type over
$S=\Spec k$ with $k$ perfect of characteristic
$p>0$,\label{thm:MKtoLRW}
\begin{equation}
C_X^{\text{M}}(n)/p^r\cong \widetilde{\nu}_{X,r}(n)(X).
\end{equation}
\end{thm}
\begin{proof} By Bloch--Kato, \cite[Theorem 2.1]{blochkato}, for any field
$F$, there is an isomorphism
$$K_n^M(F)/p^r\cong \nu^n_{F,r}(F)$$
sending $\{f_1,...,f_n\}$ to $$d\log f_1\wedge ...\wedge d\log f_n.$$
By \cite[Lemma 4.1]{grossuwa} (see
Jannsen-Saito-Sato, \cite[Theorem 2.1.1 and Theorem 2.11.3(3)]{jannsensaitosato} for
 a more detailed proof), this isomorphism respects the differentials in $C_X^{\text{M}}(n)$
and $\widetilde{\nu}_{X,r}(n)(X).$  Hence it induces an isomorphism of complexes $$C_X^\text{M}(n)/p^r \stackrel{\cong}\longrightarrow \widetilde {\nu}_{X,r}(n)(X).$$
\end{proof}

To relate $C_X^{\text{HC}}(n)$ with $C_X^\text{M}(n)$, we need the
Nesterenko--Suslin isomorphism.

\begin{thm}[Nesterenko--Suslin, \cite{NesterenkoSuslin}, Theorem 4.9]\label{isomorphism:nesterenkosuslin}Let $F$ be a field of transcendental degree $d$ over a field $k$, and $n\le d$.
Then there is an isomorphism $$\chi_F: H^{\BM}_{d+n}(F,\Z(n))\to K_{d-n}^M(F).$$
\end{thm}
\begin{proof}
The map $\chi_F$ is defined as follows: for any generator $z\in z_{n}(F,d-n)$,
$\chi_F(\bar{z})=N(\beta_z)$. Here $\bar{z}$ is the image of $z$ in
$H^{\BM}_{d+n}(F,\Z(n))$, $$N:K_{d-n}^M(k(z)) \to K_{d-n}^M(F)$$ is the norm map of Milnor
K-theory, and
$$\beta_z=\{\frac{-t_0}{t_{d-n}},..., \frac{-t_{d-n-1}}{t_{d-n}}\}\in K^M_{d-n}(k(z))$$
with $t_i$'s  the coordinates of $z$ in $\Delta^{d-n}_F$, i.e., pull-back of $t_i$ on $z$. Since $z$
intersects all the faces properly, $t_i\in k(z)^*$ for all $i$. Nestenrenko and Suslin showed that $\chi_F$ is an isomorphism.

\end{proof}

It is easy to see that this isomorphism does not depend on the choice of the base field $k$. In particular, let $F=k$, and $n\le 0$. Then $d=0$ and we have isomorphism
$$\chi_F:H^{\BM}_{n}(F,\Z(n))\cong K_{-n}^M(F).$$

\begin{lemma}\label{Lemma:HCGtoMK}
For $X$ separated and essentially of finite type over $S=\Spec k$ or
$S=\Spec D$, the Nesterenko--Suslin isomorphism induces an
isomorphism of complexes $$\chi: C_X^{\text{HC}}(n)\to
C_X^{\text{M}}(n) .$$ If the base is $\Spec k$ with characteristic $p>0$, then $$\chi/p^r:
C_X^{\text{HC}}(n)/p^r \stackrel{\cong}\longrightarrow
C_X^{\text{M}}(n)/p^r.$$ \end{lemma}

 \begin{proof} We have to show that the maps $$\chi:(C_X^{\text{HC}}(n))_{m} \to (C_X^M(n))_{m}$$ induces the following commutative diagram:
$$\xymatrix{\displaystyle \bigoplus_{X_{(m-n)}} H_{m}^{\BM}(k(x),\Z(n)) \ar[d]^-{\chi}\ar[r]^-{d''} & \displaystyle\bigoplus_{X_{(m-n-1)}}H_{m-1}^{\BM}(k(x), \Z(n)) \ar[d]^-{\chi}\\
\displaystyle \bigoplus_{X_{(m-n)}} K_{m-2n}^M(k(x)) \ar[r]^-{d'}&
\displaystyle\bigoplus_{X_{(m-n-1)}}K_{m-2n-1}^M(k(x)) .}
$$
It suffices to prove the conclusion for $m=d+n$ with $d=\dim X$, since commutativity of such diagram at the other degrees is similar. Let $X'\to X$ be the normalization of $X$, and $x$ be a codimension
1 point of $X$. Consider the following commutative diagram:
$$\xymatrix{H_{d+n}^{\BM}(k(X'),\Z(n)) \ar[r]^-{d''_{X'}} \ar@{=}[d] & \displaystyle\bigoplus_{x'|x}H_{d+n-1}^{\BM}(k(x'),\Z(n))\ar[d]^-{\sum N}\\
           H^{\BM}_{d+n}(k(X),\Z(n)) \ar[r]^-{d''_X} & H_{d+n-1}^{\BM}(k(x), \Z(n)).}$$
 Here $d''_X$ and $d''_{X'}$ are differentials in $C_X^{\text{HC}}(n)$ and $C_{X'}^{\text{HC}}(n)$, respectively, and $N=N_{k(x')/k(x)}$ is the push-forward of higher Chow groups of finite field extensions. This diagram is commutative by   covariance of Gersten complex. So $$d''_X=\sum N\circ d''_{X'}.$$

To finish the proof, consider the following diagram:
$$\xymatrix{H^{\BM}_{d+n}(k(X'),\Z(n) ) \ar[r]^-{d''_{X'}} \ar[d]^{\chi_{k(X')}} & \displaystyle\bigoplus_{x'|x} H^{\BM}_{d+n-1}(k(x'),\Z(n))\ar[d]^{\chi_{k(x')}}\ar[r] & H^{\BM}_{d+n-1}(k(x),\Z(n)) \ar[d]^{\chi_{k(x)}}\\
K_{d-n}^M(k(X')) \ar[r]^-{\partial_{x'}} & \displaystyle\bigoplus_{x'|x}
K_{d-n-1}^M(k(x'))\ar[r] & K^M_{d-n-1}(k(x)).}$$ The horizontal maps on
the right hand square are norm maps. By definition, the composition
of the  maps at the bottom is the differential $d'$ in
$C_X^{\text{M}}(n)$. By the first part of the proof, the composition
of the maps on the top is the differential $d''$ in
$C_X^{\text{HC}}(n)$. Hence it suffices to prove commutativity of
this diagram. The right hand diagram commutes by
\cite[Lemma 4.7]{NesterenkoSuslin}. On the other hand, Geisser and
Levine showed commutativity of the left hand square
\cite[Lemma 3.2]{geisserlevine} (even though their statement is
for $\Z/p$-coefficients, their proof holds for $\Z$-coefficients as
well). Hence we prove the first part. The second part about the
coefficients $\Z/p^r$ holds since both complexes are $C_X^{M}(n)$ and $C_X^{\text{HC}}(n)$ modular
$p^r$.

\end{proof}

\begin{defn}\label{definition:1} Define $$\psi=\chi\circ \phi:\Z^c_X(n)(X) \to C_X^{\text{HC}}(n)\stackrel{\cong}\to C_X^{\text{M}}(n).$$
\end{defn}

Explicitly,
\begin{defn}\label{definition:2}Denote $p:X\times \Delta^j\to X$. Given a generator $$Z\in z_n(X,m-2n)=(\Z^c_X(n)(X))_{m}, $$ we define $$\psi_m(Z)\in \bigoplus_{X_{(m-n)}}K_{m-2n}^M(k(x))=C_X^{M}(n)_m$$ (see (\ref{filtration})) as follows:

1) if $m>n+d$ or $m<n$, $\psi_m(Z)\overset{def}{=}0$.

2) if $n\le m\le n+d$ and $\dim p(Z)<m-n$,
$\psi_m(Z)\overset{def}{=}0$.

3) if $n\le m\le n+d$ and $\dim p(Z)=m-n$, then $Z$ is dominant over
some $x\in X_{(m-n)}$. Pulling back $Z$ along $$\Spec k(x)\to X,$$ we
obtain $Z_x\in z_n(k(x), m-2n)$, which is sent to $\overline{Z}_x$
by the quotient $$z_n(k(x),m-2n )\twoheadrightarrow
H_m^{\BM}(k(x),\Z(n)).$$ Applying the
Nesterenko--Suslin isomorphism $\chi_{k(x)}$, we get
$$\psi_m(Z)\stackrel{def}=\chi_{k(x)}(\overline{Z}_x)\in K_{m-2n}^M(k(x)).$$ Since
$Z$ is dominant over $x$, $Z_x$ is a closed point in
$\Delta^{m-2n}_{k(x)}$ with residue field $k(Z_x)=k(Z)$. Therefore,
by definition of $\chi_{k(x)}$, $$\psi_m(Z)=N_{k(Z)/k(x)}(\beta_Z)$$
with $t_i\in k(Z)^*$ and
$$\beta_Z=\beta_{Z_{x}}=\{\frac{-t_0}{t_{m-2n}},..., \frac{-t_{m-2n-1}}{t_{m-2n}}\}.$$
\end{defn}

 Note that in Definition \ref{definition:2}, case 3), $\overline{Z}_x=\phi_{m}(Z)$ where $\phi$
 is the map from $\Z^c_X(n)(X)$ to $C_X^{HC}(n)$.

\begin{thm}\label{Lemma:cyclecomplextoMK}
Let $X$ be separated and essentially of finite type over $\Spec k$ or
$\Spec D$. The map $\psi$ defined above is a map of complexes, and
it is functorial with respect to pullbacks defined by quasi-finite,
flat maps and push-forwards defined by proper maps.
\end{thm}
\begin{proof} Since $\psi=\chi \circ \phi$, it is a map of complexes.

For functoriality, it suffices to assume that $n\le m \le n+d$. First, we show that
$\psi$ is compatible with pull-backs defined by quasi-finite, flat
maps $f:X \to Y$. We have to prove that the following diagram is
commutative:
$$\xymatrix{z_n(Y,m-2n) \ar[d]^{f^*} \ar[r]^-{\psi_Y} &\displaystyle\bigoplus_{Y_{(m-n)}}K_{m-2n}^M(k(y)) \ar[d]^{f^{*'}} \\
 z_n(X,m-2n)\ar[r]^-{\psi_X} & \displaystyle\bigoplus_{X_{(m-n)}}K_{m-2n}^M(k(x)) .}$$
Here $f^*$ sends a generator $Z\in z_n(Y,m-2n)$ to its cycle
theoretic pull-back $$f^{-1}(Z)\in z_n(X,m-2n),$$ and
$${f^*}':K_{m-2n}^M(k(y))\to K_{m-2n}^M(k(x))$$ is defined by the
field extension $k(y)\subset k(x)$ with $x\in X_{(m-n)}$, $y\in
Y_{(m-n)}$ and $f(x)=y$. Let $$p_X:X\times \Delta^{m-2n}\to X,
p_Y:Y\times \Delta^{m-2n} \to Y$$ be the projections. If $\dim p_Y
(Z)<m-n$, then $\dim p_X (f^{-1}(Z))<m-n$, so
$$\psi_Xf^*(Z)=0=f^{*'}\psi_Y(Z).$$

Suppose that $\dim p_Y(Z)=m-n$. Without loss of generality,
replacing $Y$ by  $p_Y(Z)$ and $X$ by $X\times_{Y}p(Z)$, we can
assume that $Y$ is irreducible of dimension $m-n$ and $Z$ is
dominant over $Y$. Since $f$ is quasi-finite and flat, $X$ is of
equi-dimensional $m-n$. Let $X=\cup_iX_i$, $X_i$ be irreducible
components of $X$ and $x_i$ be the generic points of $X_i$. Then
$\dim X_i=m-n$. Therefore it suffices to prove  commutativity of the
following diagram:
$$\xymatrix{z_n(Y,m-2n) \ar[d]^{f^*}\ar[r] &   H_m^{\BM}(k(Y),\Z(n)) \ar[d]\ar[r] &  K_{m-2n}^M(k(Y)) \ar[d]^{f^{*'}}      \\
z_n(X,m-2n)\ar[r] &
\displaystyle\bigoplus_{i}H_{m}^{\BM}(k(x_i),\Z(n))\ar[r] &
\displaystyle\bigoplus_iK_{m-2n}^M(k(x_i)) }.$$ The square on the
left commutes by functoriality of higher Chow groups with respect to
flat pull back, the square on the right commutes since by definition, the
Nesterenko--Suslin isomorphism is covariant with respect to finite
field extensions. Therefore $\psi$ commutes with quasi-finite and
flat pull-backs.

If $g:X \to Y$ is a proper map, then $\psi$ is  covariant for
push-forwards. To see that, fix $$Z\in z_n(X,m-2n).$$ After replacing $X$ by the irreducible component of $X$ containing
$p_X(Z)$, it suffices to assume that $X$ is
irreducible of dimension $n$ with function field $K$ and prove that $Z$ is sent to the same element via the two paths in the following square:
$$\xymatrix{z_n(X,m-2n) \ar[d]^{g_*} \ar[r]^{\psi_X}  & K_{m-2n}^M(K) \ar[d]^{g_*'}\\
             z_n(Y,m-2n) \ar[r] ^-{\psi_Y}            & \displaystyle\bigoplus_{Y_{(m-n)}}K_{m-2n}^M(k(y)).}$$
Here $g_*(Z)$ is defined as follows:
$$g_*(Z)=\left\{\begin{array}{ll}                       0, & \text{if } \dim g(Z)<m-n;\\
                                 m_Z\cdot g(Z), & \text{if } \dim g(Z)=m-n,\end{array}\right.$$
with $m_Z=[k(Z):k(g(Z))]$, and
$$g_*'=\left\{\begin{array}{ll}    0, & \text{if } \dim g(X)<m-n;\\
                                 N_{K/k(y)}, & \text{if } X \text{ dominant over some} y\in Y_{(m-n)}. \end{array}\right.$$ To show that $g_*'\psi_X(Z)=\psi_Yg_*(Z)$, there are three cases:

1) if $\dim g(X)<m-n$, then $g_*'=0$. Moreover,  $$\dim p_Y(g_*(Z))=\dim g(p_X(Z)) \le \dim g(X)<m-n.$$
Hence $\psi_Y(g_*(Z))=0$.

2) if $X$ is dominant over some $y\in Y_{(m-n)}$ and $\dim
p_X(Z)<m-n$, then $$g_*'\psi_X(Z)=g_*'(0)=0.$$ Moreover, $\dim p_X(Z)<m-n$ also
implies $\dim p_Y(g(Z))<m-n$, hence $$\psi_Yg_*(Z)=\psi_Y(0)=0.$$

3) if $X$ is dominant over some $y\in Y_{(m-n)}$ and $\dim
p_X(Z)=m-n$, then $$\dim p_Y(g(Z))=\dim g(p_X(Z))=m-n.$$ Therefore,
$Z$ is dominant over $X$ and $g(Z)$ is irreducible and dominant over
$y$. We have a commutative diagram of field extensions:
$$\xymatrix{k(Z) & K \ar[l] \\
            k(g(Z))  \ar[u]   & k(y)\ar[u]\ar[l].}$$
Then $$\psi_Yg_*(Z)=\psi_Y(m_Z\cdot
g(Z))=N_{k(g(Z))/k(y)}(m_Z\cdot\beta_{g(Z)})$$ and
$$g_*'\psi_X(Z)=N_{K/k(y)}N_{k(Z)/K}(\beta_Z)=N_{k(g(Z))/k(y)}N_{k(Z)/k(g(Z))}(\beta_Z).$$
Since $\beta_Z$ is the image of $\beta_{g(Z)}$ under  the map
$K_{m-2n}^M(k(g(Z))) \to K_{m-2n}^M(k(Z))$,
$$N_{k(Z)/k(g(Z))}(\beta_Z)=m_Z\cdot\beta_{g(Z)}.$$ Therefore
$\psi_Yg_*(Z)=g_*'\psi_X(Z)$.

\end{proof}

\noindent\textbf{Remark 1.} ~~In \cite{landsburg}, Langsburg defined
a map from $\Z^c_X(n)(X)$ to $C_X^{\text{M}}(n)$ exactly the same as
the one in Definition \ref{definition:2}, except in case 3), instead
of using $\beta_Z$, he used
$$\beta_Z'=\{\frac{t_0}{t_{m-2n}},...,\frac{t_{m-2n-1}}{t_{m-2n}}\}.$$
By multilinearity of Milnor K-theory, it is easy to see that
$\beta_Z=\beta_Z'$ up to  2-torsion element. Therefore, the map
$\psi$ is equal to Langdsburg's map up to a 2-torsion. The advantage
of using $\beta_Z$ is that one can use the results in
\cite{NesterenkoSuslin}. On the other hand, using the idea in
Landsburg's proof of showing that his map is a map of complexes,
together with  properties of $\chi$ in \cite{NesterenkoSuslin}, we
can give a direct proof of showing that $\psi$ is a map of complexes.
This proof is lengthy, comparing to the one we give in Theorem
\ref{Lemma:cyclecomplextoMK}. (There is a small gap in Landsburg's
proof, as he only checks compatibility of his map with the
differentials in the case of discrete valuation rings).

\noindent\textbf{Remark 2.} ~~It is easy to see that the map $\psi$
can be generalized to define a map from Bloch's cycle complex $\Z^c_X(n)(X)$ to the corresponding cycle
complex with coefficients in Milnor K-groups or Quillen K-groups
defined by M. Rost \cite{rost}.

\begin{thm}\label{thm:cycletoLRW}
For any $X$ separated and essentially of finite type over $k$ of
characteristic $p>0$, the map $\psi$ induces a quasi-isomorphism
$$\hat{\psi}:\Z^c_X/p^r(n)\to \widetilde{\nu}_{X,r}(n).$$
\end{thm}

\begin{proof} Composing $\psi/p^r$ with the isomorphism
$$C_X^{\text{M}}(n)/p^r\overset{\sim}\longrightarrow
\widetilde{\nu}_{X,r}(n)(X),$$ we get a map of complexes
$$\hat{\psi}:\Z^c_X/p^r(n)(X)\to \widetilde{\nu}_{X,r}(n)(X).$$ To
compare the cohomology of $\Z^c_X/p^r(n)$ and
$\widetilde{\nu}_{X,r}(n)$, consider the niveau spectral sequence of
higher Chow groups:
\begin{equation} E^1_{s,t}=\displaystyle\bigoplus_{X_{(s)}}H_{s+t}^{\BM}(k(x),\Z/p^r(n)) \Rightarrow H_{s+t}^{\BM}(X,\Z/p^r(n)).\label{definitionspectralsequence}
\end{equation}
By \cite[Theorem 1.1]{geisserlevine}, this spectral sequence
collapses to give edge isomorphisms
\begin{equation}\Gamma: H_s^{\BM}(X,\Z/p^r(n))\cong H_s(C_X^{\text{HC}}(n)/p^r)\label{isomorphism}.
\end{equation}
Composing with the isomorphisms $$C_X^{\text{HC}}(n)/p^r \cong
C_X^{\text{M}}(n)/p^r \text{~(Lemma \ref{Lemma:HCGtoMK})}$$ and
$$C_X^{\text{M}}(n)/p^r\cong \widetilde{\nu}_{X,r}(n)(X) ~~(\text{Lemma \ref{thm:MKtoLRW}}),$$ we get an
isomorphism $$\hat{\Gamma}:H_s^{\BM}(X,\Z/p^r(n))\cong
H_s(\widetilde{\nu}_{X,r}(n)(X)).$$Since $\phi$ induces edge morphisms
of the spectral sequence, and $\psi=\chi\circ \phi$, we see that
$\hat{\psi}$ induces $\hat{\Gamma}$. Hence $\hat{\psi}$ is a
quasi-isomorphism.

\end{proof}

\noindent\textit{$\Q$-coefficients.}~~Let $k$ be a finite field with
characteristic $p$ and $X$ be smooth and projective over $k$. Let
$\Q^{c}_X(n)=\Z^c_X(n)\otimes \Q$. In \cite[Proposition 2.1]{geisser2}, conjecture A($n$) (part of Parshin's conjecture) is equivalent
to that, in the niveau spectral sequence
$$E^1_{s,t}=\displaystyle\bigoplus_{X_{(s)}}H_{s+t}^{\BM}(k(x), \Q(n))\Rightarrow H_{s+t}^{\BM}(X,\Q(n)),$$
$E^1_{s,t}=0$ for $t\neq n$. In other words, it is equivalent to the
existence of the following isomorphism:
$$H_s(E_{*,n}^1)\cong H_s^{\BM}(X,\Q(n)).$$
Here $E_{*,n}^1$ is the Gersten  complex of higher Chow groups with
$\Q$-coefficients:
$$\displaystyle\bigoplus_{X_{(d)}}H_{d+n}^{\BM}(k(x), \Q(n)) \to ...\to \bigoplus_{X_{(1)}}H_{1+n}^{\BM}(k(x), \Q(n)) \to \bigoplus_{X_{(0)}}H_n^{\BM}(k(x), \Q(n)).$$
Similar to the case above for the $\Z/p^r$-coefficients, we
obtain the following theorem:
\begin{thm}
For smooth and projective varieties over finite fields, conjecture
A(n) is true if and only if $\phi$ induces a quasi-isomorphism
  from $\Q^{c}_X(n)(X)$ to the Gersten complex of higher Chow groups $E^1_{*,n}$ in $\Q$-coefficients .
\end{thm}

\section{Comparison between Bloch's Complex and Spiess' Complex}

In this section, we assume that $S=\Spec k$ or $\Spec D$, and we will
 deal with the complex $\Z_X^c=\Z^c_X(0)$.  We
assume that $X$ is of finite type over $S$ and $\dim_S X=2$. Spiess's
dualizing complex  $\calK_X$ of \text{\'{e}}tale sheaves for surfaces uses the
weight-two motivic complex defined by \mbox{S. Lichtenbaum}
\cite[Definition 2.1]{lichtenbaum}. For a regular Noetherian ring
$A$, let $$W=\Spec A[T] ,Z=\Spec A[T]/T(T-1),$$ and
$B=\{b_1,...,b_m\}$ a finite set of exceptional units of $A$ (i.e.,
$b_i$ and $1-b_i$ are units). Let $$Y_B=\Spec
A[T]/\prod_{i=1}^m(T-b_i).$$ There is an exact sequence induced from
relative K-theory
$$K_3(A) \to K_2(W-Y_B,Z)\stackrel{{\phi_{A,B}}}\longrightarrow K_1'(Y_B) \xrightarrow{\omega_{A,B}}K_2(A). $$
Here $K_i'$ is the K-theory of the category of quasi-coherent sheaves. Taking  limits among all the $B$'s, we obtain an exact sequence
$$K_3(A) \to C_{2,1}(A)\stackrel{{\phi_A}}\longrightarrow C_{2,2}(A) \xrightarrow{\omega_A}K_2(A) $$
with $$C_{2,1}(A)=\displaystyle\underrightarrow{\lim}K_2(W-Y_B,Z), C_{2,2}(A)=\underrightarrow{\lim}K_1'(Y_B).$$ If $X$ is a
 regular Noetherian scheme, associating the group $C_{2,i}(A)$ to each regular affine scheme
 $U=\Spec A$ which is \text{\'{e}}tale over $X$, we obtain a presheaf on $X_{\text{\text{\'{e}}t}}$. Denote
 by $\underline{C}_{2,i}(X)$ the associated \text{\'{e}}tale sheaf. Then
 Lichtenbaum's weight-two motivic complex $\Z(2,X)$ is defined as the  complex
 $$\underline{C}_{2,1}(X) \stackrel{\phi_X}\rightarrow \underline{C}_{2,2}(X)$$
 with the terms in cohomological degree 1 and 2, respectively. If $A=k$ is a field, there is an exact sequence \cite[\S3]{lichtenbaum},
 \begin{equation}K_3(k)\to C_{2,1}(k)\stackrel{{\phi_k}}\longrightarrow C_{2,2}(k) \xrightarrow{\omega_k}K_2(k) \to 0 .\label{weighttwo}\end{equation}
 \begin{defn} [Spiess, \mbox{\cite[Definition 1.2.1]{spiess}}]
Define the complex of \text{\'{e}}tale sheaves $\calK_X$ as follows:
  $$\displaystyle\bigoplus_{X_{(2)}}i_{x*}\underline{C}_{2,1}(k(x)) \stackrel{c_3}\rightarrow \displaystyle\bigoplus_{X_{(2)}}i_{x*}\underline{C}_{2,2}(k(x))
  \stackrel{c_2}\rightarrow \displaystyle\bigoplus_{X_{(1)}}i_{x*}\G_m \xrightarrow{c_1}
  \displaystyle\bigoplus_{X_{(0)}}i_{x*}\Z.$$
  Here $c_3= \phi_{k(x)}$, $c_1$ is the map sending a rational function
  to its associated divisors, and $c_2$ is the composition
  $$\displaystyle\bigoplus_{X_{(2)}}i_{x*}\underline{C}_{2,2}(k(x)) \stackrel{ \omega_{k(x)}} \longrightarrow
  \displaystyle\bigoplus_{X_{(2)}}i_{x*}K_2(k(x)) \stackrel{\partial'}\rightarrow \displaystyle\bigoplus_{X_{(1)}}i_{x*}\G_m,$$ with
  $\partial'$ the map in the  Gersten resolution of K-theory \cite[\S7]{quillen}. The terms are put in (cohomological) degree -3, -2, -1, 0, respectively.\label{definition:spiesscomplex}
 \end{defn}
 Note that the assumption on degrees are different from that in \cite{spiess}.
 \begin{thm}[Spiess, \mbox{\cite[Theorem 2.2.2]{spiess}}]
 For $X$ an equidimensional surface over $\Z$ satisfying the (NR) condition and every constructible sheaf $\calF$ on $X$, there
 are perfect pairings of finite groups
 $$H^i_c(X_{\textrm{\text{\'{e}}t}},\calF)\times \Ext_X^{2-i}(\calF,\calK_X) \to\Q/\Z.$$
 \end{thm}
Here we say that $X$ satisfies the (NR) condition if
\begin{center}
\text{(NR):} ~~~~~~\textit{$k(x)$ is not formally real for every
$x\in X$.}
\end{center}
\begin{thm}[Spiess, \mbox{\cite[Proposition 2.3.2]{spiess}}]
Let $k$ be an algebraically closed field of characteristic $p$, $X$
be an irreducible surface over $k$ and $\calF$ be an $n$-torsion
constructible sheaf on $X$, with $(n,p)=1$. Let
$\calK_X(n)=RHom_X(\Z/n,\calK_X)$.  Then there are perfect pairings
of finitely generated groups:
$$H_c^i(X_{\text{\text{\'{e}}t}},\calF)\times \Ext^{1-i}_X(\calF,\calK_X(n)) \to  \Z/n.$$
\end{thm}

In \cite{deninger}, Deninger considered a dualizing complex of
\text{\'{e}}tale sheaves
\begin{equation}\calG_Y: 0 \to \displaystyle\bigoplus_{Y_{(1)}}i_{y*}\G_m \to \bigoplus_{Y_{(0)}}i_{y,*} \Z \to 0\label{deningercomplex}
\end{equation}
for curve $Y$, and E. Nart \cite{nart} compared it with  cycle complex
$\Z^c_Y$ by constructing a map from $\Z^c_Y$ to $\calG_Y$ which
induces a quasi-isomorphism. Here we put the terms in
(cohomological) degree -1 and 0 respectively. In what follows we
generalize this method and define a similar complex for surfaces.

Consider the niveau spectral sequence of higher Chow groups for
surfaces
$$E^1_{s,t}=\displaystyle\bigoplus_{X_{(s)}}H_{s+t}^{\BM}(k(x),\Z) \Rightarrow H_{s+t}^{\BM}(X,\Z).$$
Only $E^1_{0,0}$,$E^1_{1,0}$ and $E^1_{2,t}$($t\geq 0$) are
non-vanishing, $$\displaystyle H_{s+2}^{\BM}(X,\Z)\cong
\bigoplus_{X_{(2)}}H_{s+2}^{\BM}(k(x),\Z)$$ for $s> 0$, and in the
bottom of the spectral sequence, the only nonvanishing terms in
$E^1_{*,0}$ are:
\begin{equation}
\label{section4:bottomofspectralsequence}
\bigoplus_{X_{(2)}}H_2^{\BM}(k(x),\Z) \xrightarrow{f}
\bigoplus_{X_{(1)}}k(x)^* \xrightarrow{d_1}
\displaystyle\bigoplus_{X_{(0)}}\Z
\end{equation}
with
$$\coker d_1\cong H_0^{\BM}(X,\Z), \frac{\ker d_1}{\im f}\cong H_1^{\BM}(X,\Z), \ker f\cong H_2^{\BM}(X,\Z).
$$

\begin{defn}
  We define a cochain complex of abelian groups $\calC(X)$:
$$\displaystyle\bigoplus_{X_{(2)}} \frac{z_0(k(x),3)}{I_x} \xrightarrow{d_3} \displaystyle\bigoplus_{X_{(2)}}z_0(k(x),2)
\xrightarrow{d_2} \displaystyle\bigoplus_{X_{(1)}}k(x)^*
\xrightarrow{d_1} \displaystyle\bigoplus_{X_{(0)}} \Z.$$ Here
$I_x=\im(d_x:z_0(k(x),4) \to z_0(k(x),3))$, $d_3$ is the map induced
by $g_x:z_0(k(x),3) \to z_0(k(x), 2)$, $d_1=c_1$, and $d_2$ is the
composition
\begin{equation}\xymatrix{\displaystyle\bigoplus_{X_{(2)}}z_0(k(x),2) \ar@{->>}[r]^\pi \ar[dr]_{d_2} & \displaystyle\bigoplus_{X_{(2)}}H_2^{\BM}(k(x),\Z) \ar[d]^f \\
& \displaystyle\bigoplus_{X_{(1)}}k(x)^*.}
\label{definition:d_2}\end{equation} The (cohomological) degrees of the terms are -3, -2, -1, 0, respectively.
\end{defn}
\begin{thm}
Let $X$ be a surface over $\Spec k$ or $\Spec D$, then there is a
map $$\psi':\Z^c_X(X) \to \calC(X)$$ which induces isomorphisms
$$H^{-i}(\Z^c_X(X))\cong H^{-i}(\calC(X))$$
for $0\leq i\leq 3$\label{theoremcycletruncated}.
\end{thm}
\begin{proof} Following the same idea of the definition of $\psi$ in \S2,
we can define a map from $\tau_{\ge -3}\Z^c_X(X)$ to $\calC(X)$:
$$\xymatrix{ z_0(X,3)/I' \ar[r]^-d\ar[d]^-{\psi_3'} & z_0(X,2) \ar[r]^-d\ar[d]^-{\psi_2'} & z_0(X, 1) \ar[d]^-{\psi_1'}\ar[r]^-d & z_0(X,0) \ar[d]^-{\psi_0'}.
\\
\displaystyle\bigoplus_{X_{(2)}} z_0(k(x),3)/I_x \ar[r]^-{d_3} &
\displaystyle\bigoplus_{X_{(2)}}z_0(k(x),2) \ar[r]^-{d_2} &
\displaystyle\bigoplus_{X_{(1)}}k(x)^* \ar[r]^-{d_1} &
\displaystyle\bigoplus_{X_{(0)}} \Z}$$ Here $$I'=\im(z_0(X,4) \to
z_0(X,3)), \psi_0'\stackrel{def}=\psi_0, \psi_1'\stackrel{def}=\psi_1.$$ For $i=2,3$, $\psi_i'$ is defined as
follows: if $\dim p(Z)\le 1$, $\psi_i'(Z)=0$; if $\dim p(Z)=2$, then
$Z$ is dominant over some $x\in X_{(2)}$, so pulling back $Z$ along
$\Spec k(x)\to X$, we get an element $\psi_i'(Z)\in z_0(k(x),i)$. By
similar argument as in \S2, we see that the diagram is commutative.
Moreover, the $\psi_i'$'s induce the corresponding isomorphisms from the
degeneration of the niveau spectral sequence: $$H_0^{\BM}(X,\Z)\cong
\coker d_1\cong H_0(\calC(X)),$$
$$H_1^{\BM}(X,\Z) \cong \frac{\ker d_1}{\im f}\stackrel{(*)}=\frac{\ker d_1}{\im d_2}= H_1(\calC(X)),$$
$$H_3^{\BM}(X,\Z) \cong \bigoplus_{X_{(2)}}H_3^{\BM}(k(x),\Z)\cong \ker d_3\cong H_3(\calC(X)),$$
$$H_2^{\BM}(X,\Z) \cong \ker f .$$
Here the equality $(*)$ follows from the diagram (\ref{definition:d_2}).

It remains to show that $\ker f\cong H_2(\calC(X))$. Look at the diagram (\ref{definition:d_2}). There is an exact
sequence
$$0 \to \ker \pi \to \ker d_2 \to \ker f \to \coker \pi \to \coker d_2\to \coker f \to 0.$$
Since $\pi$ is surjective, $\coker \pi=0$, so $\ker f\cong \ker
d_2/\ker \pi$. But $\ker \pi=\im d_3$. Hence
$$H_2^{\BM}(X,\Z)\cong \ker f\cong \ker d_2/\im d_3\cong H_2(\calC(X)).$$

\end{proof}

\noindent\textbf{Remark.}~ ~In the above proof, $\psi_1$ is similar to
the map defined by Nart (\cite{nart}). As noted by Nart, there are
only two types of generators in $z_0(X,1)$, the vertical ones and
the horizontal ones. The vertical ones are those $Z$'s with $\dim
p(Z)=0$, and the horizontal ones are those $Z$'s with
$p(Z)=\overline{\{x\}}$ for some $x\in X_{(1)}$. Nart's map sends
the first type to 0 and the second type to $N(\frac{t_0}{t_1})$ with
$N:k(Z)^*\to k(x)^*$ the norm map of fields. In our situation,
according to the Nesterenko--Suslin isomorphism, $\psi_1$ sends the
first type to 0 and the second type to $N(\frac{-t_0}{t_1})$. So it agrees with Nart's map up to  2-torsion.

 \begin{thm}\label{thm:calCcalK}
   For a surface $X$ over $\Spec k$ or $\Spec D$, $\calC(X)\cong \calK_X(X)$ in the derived category of \text{\'{e}}tale sheaves.
 \end{thm}

 Before proving the theorem, we need the following lemma. First, let us recall some notations from
\cite{blochlichtenbaum}.  From now on until Theorem \ref{thmcyclespiess},
assume that $F$ is a field with $\dim_SF=2$. Let $\Delta^p=\Delta^p_S$. A closed subvariety
$$\sigma:t_{i1}=t_{i2}=...=t_{iq}=0$$ is called a codimension $q$
face. A closed subvariety $V\subset \Delta^p$ is said to be in good
position if $V\cap \sigma$ has codimension $\ge q$ in $V$ for any
codimension $q$ face $\sigma$. Let $$\calV^n=\calV^n(\Delta^p)\subset
\Delta^p$$ denote the union of all codimension $n$ closed
subvarieties of $\Delta^p$ in good position. Given a scheme $X$, we
write $K(X)$ for some space, contravariant functorial in $X$, whose homotopy
groups calculate the Quillen K-theory of $X$. For $Y\subset X$ a
closed subset, we write $K(X,Y)$ for the homotopy fibre of $K(X) \to
K(Y)$. This construction can be iterated. Given $Y_1,...,Y_n\subset
X$, we define the multi-relative K-space
$$K(X;Y_1,...,Y_n)\stackrel{def}=\text{homotopy fibre of}$$
$$\Big(K(X;Y_1,...,Y_{n-1}) \to K(Y_n;Y_1\cap Y_n,...,Y_{n-1}\cap Y_n)\Big).$$
Let $$K(\Delta^p,\partial)=K(\Delta^p;\partial_0,...,\partial_p)$$
$$K(\Delta^p,\sum)=K(\Delta^p;\partial_0,...,\partial_{p-1}),$$
where $\partial_i=\partial_i(\Delta^{p-1})$. Let $\Psi=\partial$ or
$\sum$. Define
$$K_V(\Delta^p,\Psi)\stackrel{def}=\text{homotopy fibre of }(K(\Delta^p,\Psi) \to K(\Delta^p-V,\Psi-V)).$$
If $W\subset V$ is an inclusion of closed subvarieties, then there
is a canonical map
$$K_W(\Delta^p,\Psi)\to K_V(\Delta^p,\Psi). $$
Hence we can define
$$K_{\calV^n}(\Delta^p,\Psi)\stackrel{def}=\underset{V\subset \calV^n}{\underrightarrow{\lim}}K_V(\Delta^p,\Psi)$$
$$K_{\calV^n-\calV^{n+1}}(\Delta^p,\Psi)\stackrel{def}=\underset{V\subset \calV^n}{\underrightarrow{\lim}} \underset{W\subset \calV^{n+1}}{\underrightarrow{\lim}}K_{V-W}(\Delta^p -W, \Psi-W).$$

 There are two isomorphisms from $H_{2}^{\BM}(F,\Z)$ to $K_2(F)$ in \cite{blochlichtenbaum}: the first one $\theta'$ is induced by the edge morphisms of
the spectral sequence
\begin{equation}E_2^{p,q}=H_{q-p+4}^{\BM}(F,\Z(2+q))\Rightarrow K_{-p-q}(F)\label{spectralsequence:motivic}
 \end{equation}
 proved in \cite{blochlichtenbaum} (note that we have changed the notation in loc. it. to the motivic Borel--Moore homology of $F$ over base $S$); the second one $\theta$ is induced from the quasi-isomorphism between $T_F$ and $\Z(2,F)$ (loc.it., \S7). Here $T_F$ is the following truncation of cycle complex of fields:
 $$z_0(F,3)/I_F \to z_0(F,2)$$
 with $$I_F\stackrel{def}=\im(z_0(F,4)\to z_0(F,3)).$$ Hence $\coker T_F=H_{2}^{BM}(F,\Z)$.  From (\ref{weighttwo}), we know that
 this quasi-isomorphism induces an isomorphism:
\begin{equation}\theta:H_{2}^{BM}(F,\Z)\to K_2(F).\label{theta}\end{equation}
\begin{lemma}\label{lemma:BL}
 (1) For any field $F$, $$\theta=\theta':H_{2}^{BM}(F,\Z)\to K_2(F).$$

 (2) The following diagram is commutative
 $$\xymatrix{H_2^{\BM}(K,\Z) \ar[rd]_{\chi} \ar[d]^{\theta}& \\
K_2(K) \ar[r]_{s} & K_2^M(K) .}$$
\end{lemma}
\begin{proof} (1)$\Rightarrow$(2): By
\cite[Proposition 3.3]{geisserlevine}, we know that $s\circ
\chi=\theta'$. Then from (1), we have $$s\circ \chi=\theta'=\theta.$$

(1) By the definitions of $\theta$ and $\theta'$ \cite[\S1 and \S7]{blochlichtenbaum}, and the following isomorphisms:
$$K_0(\Delta^2,\partial)\stackrel{\cong}\longleftarrow K_1(\Delta^1,\partial)\stackrel{\cong}\longleftarrow K_2(\Delta^0), $$
we see that the two maps are induced by the two paths from
$H_2^{\BM}(F,\Z)$ to $K_1(\Delta^1,\partial)$ and to
$K_0(\Delta^2,\partial)$ shown in the following diagram:
$$\xymatrix{
 K_{1,\calV^1-\calV^2}(\Delta^2,\sum)\ar@{->>}[r] \ar@{->>}[d] & K_{0,\calV^2}(\Delta^2,\sum) \ar[d]^{\cong} &  K_{0,\calV^2}(\Delta^2,\partial) \ar[l]^{\cong} \ar@{->>}[d]\\
 K_{1,\calV^1}(\Delta^1,\partial) \ar@{->>}[d] & z_0(F,2) \ar@{->>}[d] &  K_{0,\calV^1}(\Delta^2,\partial)\ar[d]^{\cong} \\
 K_{1}(\Delta^1,\partial)& H_{2}^{BM}(F,\Z) & K_0(\Delta^2,\partial) .}$$
 The map $\theta'$ is the path from $z_0(F,2)$ to $K_{0,\calV^2}(\Delta^2,\partial)$ and then to $K_0(\Delta^2,\partial)$, while $\theta$ is the path from $z_0(F,2)$ to $K_1(\Delta^1,\partial)$ via $K_{1,\calV^1-\calV^2}(\Delta^2,\sum)$. In order to show that $\theta=\theta'$,
it suffices to show commutativity of the following diagram:
(starting from $K_{1,\calV^1-\calV^2}(\Delta^2,\sum)$ and end at
$K_0(\Delta^2,\partial)$)
$$\xymatrix{K_{0,\calV^2}(\Delta^2,\sum) & K_{0,\calV^2}(\Delta^2,\partial)\ar[l]^-{\cong}_-3 \ar@{->>}[r]^{4}  & K_{0,\calV^1}(\Delta^2,\partial) \ar@{->>}[r]& K_0(\Delta^2,\partial) \\
K_{1,\calV^1-\calV^2}(\Delta^2,\sum) \ar@{->>}[rr]^2\ar@{->>}[u]^1 &
& K_{1,\calV^1}(\Delta^1,\partial) \ar@{->>}[u]^5 \ar@{->>}[r]&
K_1(\Delta^1,
\partial)\ar[u]^{\cong}}$$
Here the right hand square is induced from the embedding
$(\Delta^1,\partial_0,\partial_1)\subset
(\Delta^2,\partial_0,\partial_1)$, so it is commutative. Therefore
we have to show that the square on the left commutes, or more
precisely, $5\circ 2=4 \circ 3^{-1}\circ 1.$

Recall the long exact sequence of relative K-groups:
\begin{equation}\label{sequence:relativektheory} K_1(Y) \stackrel{i}\to K_0(X,Y) \stackrel{j}\to K_0(X) \stackrel{k}\to K_0(Y).
\end{equation}
That is the only type of long exact sequence involved in the
diagrams above. In particular,  $$1\sim i, 2\sim k, 3\sim j, 4\sim
k, 5\sim i. $$ Here  map $1$ is for the embedding
$(\Delta^2-\calV^1,\sum)\subset (\Delta^2-\calV^2,\sum)$, map 4 is for
the embedding $(\Delta^2-\calV^1,\partial)\subset (\Delta^2-
\calV^2,\partial)$,  maps 2, 3 and 5 are for the embedding
$(\Delta^1,\partial)\subset (\Delta^2,\sum)$ (embedding of the last
face), and $1\sim i$ means that the map $1$ corresponds to type $i$ in the sequence (\ref{sequence:relativektheory}), etc. By functoriality of relative sequence of K-theory, we obtain
the commutativity.
\end{proof}

\begin{proof} [Proof of Theorem \ref{thm:calCcalK}] Let us first define a map $\calC(X)\to \calK_X(X)$ in the derived category.
 For any $x\in X_{(2)}$, there is a quasi-isomorphism \cite[Theorem 7.2]{blochlichtenbaum} in the derived category between the complex $$z_0(k(x),3)/I_x \to z_0(k(x), 2)$$ and
  the complex $$C_{2,1}(k(x))\to C_{2,2}(k(x)).$$
  Moreover, it induces an isomorphism between  the cokernels of the two complexes $$H_2^{\BM}(k(x),\Z)\overset{\theta_x}{\underset{\sim}{\longrightarrow}} K_2(k(x))$$(see (\ref{theta}) above).
  Recall that the complex $\calC(X)$ (resp., $\calK_X(X)$) is defined by connecting
  $$\displaystyle\bigoplus_{X_{(2)}} (z_0(k(x),3)/I_x) \xrightarrow{d_3} \displaystyle\bigoplus_{X_{(2)}}z_0(k(x),2)$$
  $$\text{\Big(resp., }\displaystyle\bigoplus_{X_{(2)}}C_{2,1}(k(x)) \stackrel{c_3}\rightarrow \displaystyle\bigoplus_{X_{(2)}}C_{2,2}(k(x)) \text{\Big)}$$
  with
  $$\displaystyle\bigoplus_{X_{(1)}}k(x)^* \xrightarrow{d_1} \displaystyle\bigoplus_{X_{(0)}} \Z.$$
  via $\oplus_{X_{(2)}}H_2^{\BM}(k(x),\Z)$ (resp., $\oplus_{X_{(2)}}K_2(k(x))$). Since $c_1=d_1$, it suffices to show that the following diagram is (anti-)commutative:
    $$\xymatrix{
\displaystyle\bigoplus_{X_{(2)}}H_2^{\BM}(k(x), \Z)\ar[r]^-f \ar[d]^-{\theta_x}  &  \displaystyle\bigoplus_{X_{(1)}}k(x)^* \ar@{=}[d] \\
\displaystyle\bigoplus_{X_{(2)}}K_2(k(x))\ar[r]^{\partial'} &
\displaystyle\bigoplus_{X_{(1)}}k(x)^* .}$$
 Let $\eta\in X_{(2)}$ and $K=k(\eta)$,
  consider the following diagram:
$$\xymatrix{H_2^{\BM}(K,\Z) \ar[rd]_{\chi}\ar[drr]^f \ar[dd]^{\theta}  & & \\
& K_2^M(K)\ar[r]^-{\partial}\ar[dl]_s  & \displaystyle\bigoplus_{X_{(1)}}k(x)^* \\
K_2(K) \ar[rru]_{\partial'} & &  .}$$ Here $\partial$ is the tame symbol and $\chi$ is the
Nesterenko--Suslin isomorphism. We have to show that the outside
triangle is (anti-)commutative, so it suffices to show that the
three small triangles are (anti-)commutative. The triangle on the
top is commutative by Lemma \ref{Lemma:HCGtoMK}. The one on the left commutes by Lemma \ref{lemma:BL}.

For the triangle on the bottom, we have to show commutativity of the
following diagram:
$$\xymatrix{K_2^M(K) \ar[r]^-\partial \ar[d]^s & k^*\\
K_2(K)\ar[ru]_{\partial'} &  .}$$ Here $k$ is the residue field of a
valuation $v$ of $K$, with $\pi$ a prime element. From \cite[Definition 1.1 R3e]{rost}, $\partial'$ has the following property: for any
$\rho, u\in K^*$  with $v(u)=0$, $\partial'(\{u\}\cdot
\rho)=-\{\bar{u}\}\cdot v(\rho)$, i.e. $\partial'\circ
s(\{u,\rho\})=\{\bar{u}\}^{-v(\rho)}$. Here $\bar{u}$ is the residue
class of $u$ in $k^*$. But
$$\partial(\{u,\rho\})=(-1)^{v(u)v(\rho)}\overline{\{\frac{u^{v(\rho)}}{\rho^{v(u)}}\}}=\{\bar{u}\}^{v(\rho)}.$$
By multilinearity of Milnor K-theory,  any element
$\{u_1\pi^n,u_2\pi^m\}$ with $v(u_i)=0$ can be decomposed into a
product of elements of the form $\{u,\pi\}$ and
$\{\pi,\pi\}=\{-1,\pi\}$ with $v(u)=0$. Therefore,  the diagram
commutes up to sign. So the triangle outside is commutative. In
conclusion, there is an morphism from $\calC(X)$ to $\calK_X$ in the
derived category of \text{\'{e}}tale sheaves.

  Now let us compare the (co)homology groups.
  It is clear that $$H_0(\calC(X))\cong H_0(\calK_X(X)).$$ From \cite[Theorem 7.2]{blochlichtenbaum},  $$H_3(\calC(X))\cong H_3(\calK_X(X)).$$  To compare $c_2$ and $d_2$, consider the diagrams defining them
  $$\xymatrix{\displaystyle\bigoplus_{X_{(2)}}z_0(k(x),2) \ar[r]^\pi \ar[dr]_{d_2} & \displaystyle\bigoplus_{X_{(2)}}H_2^{\BM}(k(x),\Z) \ar[d]^f  & \displaystyle\bigoplus_{X_{(2)}}C_{2,2}(k(x)) \ar[r]^{ \omega_{k(x)}}\ar[dr]_{c_2}& \displaystyle\bigoplus_{X_{(2)}}K_2(k(x))\ar[d]^{\partial'}\\
& \displaystyle\bigoplus_{X_{(1)}}k(x)^* &
&\displaystyle\bigoplus_{X_{(1)}}k(x)^* .}$$ Here $\pi$ and $
\omega_{k(x)}$ are both surjective onto groups that are connected via the isomorphism
$$\theta:H_2^{\BM}(k(x), \Z) \to K_2(k(x)),$$
and $\partial'\theta=f$, so $$\im d_2=\im f=\im\partial' =\im c_2.$$
Therefore,
$$H_1(\calC(X))\cong \ker d_1/\im d_2\cong \ker c_1/\im c_2 \cong H_1(\calK_X(X)).$$
Similar to the proof of Theorem \ref{theoremcycletruncated}, we have
that $$H_2(\calC(X))\cong \ker f\cong \ker \partial' \cong
H_2(\calK_X(X)) . $$

\end{proof}

In the following theorem we will need to consider fields of dimension $2$ over the base $S$.
If
$S=\Spec k$, a field $K$ of transcendental degree 2 over $S$ has $\dim_SK=2$. If $S=\Spec D$,
the function field $K$ of a scheme $X$ with $\dim_SX=2$ will have $\dim_SK=2$.

\begin{thm}\label{thmcyclespiess}
Suppose that $X$ is a two-dimensional scheme over $S=\Spec k$ or $S=\Spec
D$ and the function field of $X$ is $K$. Then there is a map in the derived category $\Z^c_X\to \calK_X$ which induces a quasi-isomorphism
 $$\tau_{\ge-3}\Z^c_X \cong \calK_X.$$
 Assuming the Beilinson--Soul\text{\'{e}} Conjecture for $\Z(2)$ on fields of dimension $2$ over $S$, then $\Z^c_X\cong \calK_X$ in the derived category. Here $\Z(n)$ is Voevodsky's motivic complex.
\end{thm}

\begin{proof}   From Theorem \ref{theoremcycletruncated} and
\ref{thm:calCcalK}, we know that in the derived category, there is a
map $\Z^c_X \to  \calK_X$.
Moreover, $$\tau_{\ge -3}\Z^c_X\cong \calK_X.$$

The Beilinson--Soul\text{\'{e}} Conjecture \cite[Conjecture 5]{kahn} asserts that for a regular scheme $X$ and $n>0$, $H^i(X,\Z(n))=0$
for $i\le 0$. In particular, for $\Z(2)$ on the function field $K$ of $X$ (note that $\dim_SK=2$), it is equivalent to asserting that $$H_i^{\BM}(K,\Z)=H^{4-i}(K,\Z(2))=0$$ for $i\ge 4$. To prove that $\Z^c_X\cong \tau_{\ge -3}\Z^c_X$, we need to show that
$\Z_X^c(X)$ is acyclic at (cohomological) degree $\le -4$. Now we have two cases,
\textit{(i)} $X$ over $S=\Spec k$, or $X$ over $S=\Spec D$ but factors through some closed point,
\textit{(ii)} $X$ flat over $S=\Spec D$.

\textit{case (i)} This means that $X$ is a scheme over a perfect field $k$. Let $U \hookrightarrow X$ be any open subscheme of $X$ with complement
$Z\hookrightarrow X$. Enlarge $Z$ so that it has dimension $1$. Then there is a long exact sequence
$$\to H_i^{\BM}(Z,\Z) \to H_i^{\BM}(X,\Z) \to H_i^{\BM}(U,\Z) \to.$$
By \cite[Theorem 1.1]{nart}, $\Z^c_Z(Z)$ is quasi-isomorphic to $\calG_Z(Z)$ (Deninger's complex,
see (\ref{deningercomplex})), so $H_i^{\BM}(Z,\Z)=0$ for $i\ge 2$. Hence $H_i^{\BM}(X,\Z)\cong H_i^{\BM}(U,\Z)$ for $i\ge 3$. Take limit among all the open subschemes over $X$, we get that, for $i\ge 3$, $$H^{-i}(\Z_X^c(X))=H_i^{\BM}(X,\Z)\cong H_i^{\BM}(K, \Z)=0.$$ Hence $$\Z^c_X\cong \tau_{\ge -3}\Z^c_X\cong \calK_X.$$

\textit{case (ii)} $X$ is flat over $S=\Spec D$. This is a local problem, so we can assume that $D$ is a discrete valuation ring. Let $U$ be the generic fiber and $Z$ be the closed fiber. Then $Z$ is a curve over the residue field of $D$. Similar as in \textit{(i)}, we can get the conclusion.
\end{proof}

\begin{cor}
\label{corollaryRhom}
Suppose that $X$ is a two-dimensional scheme over $S=\Spec k$ or $S=\Spec
D$. Then for any torsion sheaf $\calF$ on $X$, there is an isomorphism
  $$\RHom_X(\calF,\Z^c_X)\to \RHom_X(\calF,\calK_X).$$
\end{cor}
\begin{proof} First, we assume that the conclusion is true for constructible
sheaves on surfaces smooth over $S=\Spec k$ or $\Spec D$. If $X$ is over $S=\Spec k$ and not smooth, and $\calF$ is constructible, choose an open
smooth subscheme $U\stackrel{j}\to X$ with complement
$Z\stackrel{i}\to X$ of codimension 1. By purity, $Ri^!\Z_X^c\cong
\Z^c_Z$ \cite[Proposition 3.5]{geisser} and $Ri^!\calK_X\cong
\calG_Z$ \cite[Proposition 1.3.1 (ii)]{spiess}.  From the short
exact sequence $$0 \to j_!j^*\calF \to \calF \to i_*i^* \calF \to 0,$$
we obtain a commutative diagram of distinguished triangles
\begin{equation}\xymatrix{\RHom_X(i_*i^*\calF, \Z_X^c)\ar[r]\ar[d] & \RHom_X(\calF, \Z^c_X) \ar[d]^{\epsilon_X}\ar[r] & \RHom_X(j_!j^*\calF, \Z_X^c)\ar[d]\\
            \RHom_X(i_*i^*\calF,\calK_X) \ar[r]    & \RHom_X(\calF,\calK_X) \ar[r] & \RHom_X(j_!j^*\calF, \calK_X).}
\label{diagramlocalization}
\end{equation}
By adjointness and purity, the groups on the left are canonically
isomorphic to $$\RHom_Z(i^*\calF,\Z_Z^{c}) \text{~and~}
\RHom_Z(i^*\calF,\calG_Z),$$ respectively. Hence, by Nart
\cite{nart}, the left hand vertical map is an isomorphism. By
adjointness, the groups on the right  can be identified with
$$\RHom_U(j^*\calF,\Z_U^{c})\text{~and~} \RHom_U(j^*\calF,\calK_U),$$
respectively. Since $U$ is smooth, the right hand vertical map is an
isomorphism by induction assumption. Therefore, the vertical map in the middle
is an isomorphism.

Let $X$ be a scheme over $S=\Spec D$ but not smooth and $\calF$ is constructible. If $X\to S$ factors through some
closed point of $S$, then it is indeed a variety over a perfect field, in which case it is proved as above. If $X\to S$ does not factor through any closed point of $S$, then it is flat over $S$. It suffices to assume that $D$ is a discrete valuation ring. Let $U$ be the generic fiber and
$Z$ be the special fiber, then $Z$ is a curve over a field. For the decomposition $$U\subset X\supset Z,$$ we can apply the strategy as above to prove the conclusion for $X$ (see \cite[Corollary 7.2]{geisser} for purity of cycle complex).

For general torsion sheaf $\calF$, it can be written as
$\calF=\colim \calF_i $ with $\calF_i$ constructible \cite[Chapter II, Proposition 0.9]{milne1}. So $$\RHom_X(\calF_i,\Z^c_X)\cong
\RHom_X(\calF_i,\calK_X).$$ Since for any sheaf $\calH$, $$\RHom_X(\colim\calF_i,\calH)\cong
\text{R}\lim\RHom_X(\calF_i,\calH),$$   the
canonical map $\Z^c_X\to \calK_X$ in the derived category induces a
map
$$\text{R}\lim\RHom_X(\calF_i,\Z^c_X) \to \text{R}\lim\RHom_X(\calF_i,\calK_X),$$
hence a map of spectral sequences
$$\xymatrix{E_2^{s,t}=\lim^s\Ext^t_X(\calF_i,\Z^c_X) \ar@{=>}[r]\ar[d] & \Ext^{s+t}_X(\calF,\Z^c_X)\ar[d]^{\epsilon_X}\\
            E_2^{s,t}=\lim^s\Ext^t_X(\calF_i,\calK_X) \ar@{=>}[r]       & \Ext^{s+t}_X(\calF,\calK_X).}$$
Note that these two spectral sequences are convergent. For the first
one, it follows from \cite[Corollary 4.9]{geisser} (the statement in
this Corollary is for $\Spec k$-schemes, but the proof holds for
$\Spec D$-schemes). For the second one, since $\calK_X\cong
\tau_{\ge -3}\Z^c_X$, hence it is a bounded spectral sequence. The two spectral sequences are isomorphic at
the $E_2$-page via maps defined by $\Z^c_X\to \calK_X$, hence
$$\RHom_X(\calF,\Z_X^c)\cong \RHom_X(\calF,\calK_X).$$

Now let us prove the conclusion for constructible sheaves on
surfaces smooth over $S$.  Again, we separate the cases: a) $S=\Spec k$ with
characteristic $p$; b)  $S=\Spec D$.

 a) Let $\calF$ be a $\Z/p^r$-sheaf, with $p=\text{char}(k)$. By the niveau spectral sequence of higher Chow groups with $\Z/p^r$ coefficients, for $t\ge 3$, $$H_{t}^{\BM}(X,\Z/p^r)\cong \displaystyle\bigoplus_{X_{(2)}}H_t^{\BM}(k(x),\Z/p^r).$$ By \cite[Theorem 1.1]{geisserlevine},
we see that $H_{t}^{\BM}(X,\Z/p^r)=0$ for $t\ge 3$. Therefore,
$$\Z_X^c/p^r\cong \tau_{\ge -3} \Z_X/p^r\cong \calK_X/p^r.$$ In
particular, for such $\calF$,
  $$\epsilon_X: \RHom_X(\calF,\Z^c_X)\stackrel{\cong}\longrightarrow \RHom_X(\calF,\calK_X).$$
Now let $\calF$ be an $n$-torsion sheaf, with $(n,p)=1$. Let
$\Z/n(2)=\mu_n\otimes \mu_n$. Since both $\Z^c_X$ and $\calK_X$
satisfy the Kummer sequence (see \cite[Theorem 1.5]{geisserlevine2} and \cite[Proposition 1.5.1]{spiess}, respectively):
$$\to \Z^c_X\stackrel{n}\to \Z^c_X \to \Z/n(2) \to,$$
$$\to \calK_X\stackrel{n}\to \calK_X \to \Z/n(2) \to,$$
thus, $$\RHom_X(\calF,\Z^c_X)\cong \RHom_X(\calF, \calK_X).$$

b) Let $X$ be surface smooth over $\Spec D$.  Let $\calF$ be an
$n$-torsion. If $n$ is invertible on $X$, then $\Z^c_X$ and
$\calK_X$ satisfy the Kummer sequence \cite[Proposition 1.5.1]{spiess}. Hence $$\epsilon_X:\RHom(\calF,\Z^c_X)\to
\RHom(\calF,\calK_X)$$ is an isomorphism. If $n$ is not invertible on $X$, consider the following
decomposition
$$\xymatrix{W\ar[r]^{i}\ar[d] & X \ar[d] & V \ar[l]_{j}\ar[d] \\
            \Spec D/(n) \ar[r] & S= \Spec D & \Spec D[1/n] .\ar[l]}$$
            Here the two squares are fiber products. Then $i$ and $j$ are closed and open immersions, respectively. Since $X$ is flat over $S$, $W$ has codimension 1. Similar as in (\ref{diagramlocalization}), there is a localization diagram and the vertical maps on the left and right are isomorphisms, which proves that
             $$\epsilon_X:\RHom_X(\calF,\Z^c_X)\to \RHom_X(\calF,\calK_X)$$ is an isomorphism. This finishes the proof.

\end{proof}

\section{Comparison between Bloch's Complex and Sato's Complex}

\noindent\textit{Sato's complex $\frakI_r(n)_X$.}~~K. Sato
\cite{sato} defined a complex $\frakI_r(n)_X$ for certain schemes
over Dedekind domains. Let us first recall his notations and
definitions. In this section, $S=\Spec D$. Let $\Sigma$ be the set of closed points of $S=\Spec D$
of characteristic $p$. Suppose that $\Sigma\neq \emptyset$. Throughout
this section, the scheme $X$ over $D$ is assumed to satisfy the
following condition:

\begin{cond}[\cite{sato}, Condition 4.1.2]\label{condition*}Let $X$ be a regular, pure-dimensional scheme of dimension $d$, flat and of finite type over $S$. For any $s\in \Sigma$, each connected component $X'$ of $$X\times_S\Spec(D^{h}_{s})$$ is a regular semi-stable family over $D^h_s$ or over the integral closure of $D^h_s$ in $\Gamma(X',\calO_{X'})$.
Here $D^h_s$ is the Hensenlization of $D$ at $s$, and a scheme flat and
of finite type over a discrete valuation ring is called a regular
semi-stable family if it is regular, its generic fiber is smooth, and its
special fiber is a normal crossing divisor.
\end{cond}

Let $Y\subset X$ be the divisor defined by the radical ideal $(p)\subset
\calO_X$. Let $U$ be its complement, and $\iota$ and $j$ be as
follows:
$$\xymatrix{U \ar[r]^{j} & X & Y \ar[l]_{\iota}.}$$

\begin{defn}
\label{definitionof nu}Define $\nu_{Y,r}^n$ to be the \text{\'{e}}tale sheaf on $Y$
$$\nu^n_{Y,r}\stackrel{def}=\ker \Big(\bigoplus_{y\in Y^{(0)}}i_{y*}\nu^n_{y,r}
\stackrel{\partial}\longrightarrow
\bigoplus_{Y^{(1)}}i_{y*}\nu^{n-1}_{y,r}\Big).$$ Here
$\partial$ is the map in the Gersten complex of logarithmic de
Rham-Witt sheaves (see \cite[\S1]{kato2} or Lemma \ref{thm:MKtoLRW} above).
\end{defn}

 Define $M_r^n=\iota^*R^nj_*\mu_{p^r}^{\otimes n}$ on $Y$. We recall Sato's definition of the map $$\sigma^n_{X,r}:M_r^n \to \nu^{n-1}_{Y,r}.$$ Since it is a map of sheaves, it can be defined  on the stalks at points in $Y$, i.e., points of characteristic
$p$. Let $y\in Y$ and $s\in \Sigma$ be its image in $S$. After replacing $D$ by $D_s^h$ and $X$ by $X\times_S\Spec(D^h_s)$ (and then replacing $D_s^h$ by its integral closure in one of the irreducible components of $X\times_S\Spec(D_s^h)$, if necessary, depending the condition (\ref{condition*})), we can assume that $X$ is a regular semi-stable family over an Henselian discrete valuation ring $O_K$ with function field $K$ of characteristic $0$ and residue field $k$ of characteristic $p$.
Let $$i_x:x\to X_K \to X,  i_y:y\to X_k$$ be the embedding of points in $X_K$ and $X_k$, respectively. Then there is a diagram:
\begin{equation}\label{definition:sigma}
\xymatrix{ & M_r^n \ar[r] & \displaystyle\bigoplus_{X_K^{(0)}}\iota^*R^n( i_{x})_*\mu_{p^r}^{\otimes n} \ar[d]^-{\partial_1}\ar[r]^-{\partial_2} & \displaystyle\bigoplus_{X_K^{(1)}}\iota^*R^{n-1}( i_{x})_*\mu^{\otimes n-1}_{p^r}\ar[d]^-{\partial_3}\\
0 \ar[r]& \nu_{X_k,r}^{n-1} \ar[r] &
\displaystyle\bigoplus_{X_k^{(0)}}i_{y*}\nu_{y,r}^{n-1}
\ar[r]^-{\partial} &
\displaystyle\bigoplus_{X_k^{(1)}}i_{y*}\nu_{y,r}^{n-2}.}
\end{equation}
Here all the $\partial_i$'s are boundary maps defined in \cite[\S1]{kato2} (or see \cite[Lemma 3.2.4]{sato}). The upper row is part of the Gersten complex of $\mu_{p^r}^{\otimes n}$ on $X$. The lower row is exact by definition. The square is anti-commutative
by Kato, \cite[1.7]{kato2}. Hence it induces a map $$\sigma_{X,r}^n:
M_r^n \to \nu^{n-1}_{X_k,r}.$$

 Define the map $\sigma_{X,r}(n)$ as the composition:
$$\sigma_{X,r}(n):\tau_{\le n}Rj_*\mu_{p^r}^{\otimes n} \to (R^nj_*\mu_{p^r}^{\otimes n})[-n]\overset{\sim}{\underset{\diamond}\longrightarrow}(\iota_*M_r^n)[-n] \stackrel{\sigma^n_{X,r}}\longrightarrow \iota_*\nu^{n-1}_{Y,r}[-n].$$
Here $\diamond$ is the adjunction $Id \to \iota_*\iota^*$.

\begin{defn}[Sato, \mbox{\cite[Lemma 4.2.2]{sato}}]\label{definition:satocomplex}
If $n\ge 1$, and let $\frakI_r(n)_X$ be a complex fitted into the
following distinguished triangle
$$\xymatrix{\iota_*\nu^{n-1}_{Y,r}[-n-1] \ar[r]^-g& \frakI_r(n)_X \ar[r]^-t & \tau_{\le n}Rj_*\mu_{p^r}^{\otimes n} \ar[rr]^{\sigma_{X,r}(n)} & & \iota_*\nu^{n-1}_{Y,r}[-n].}$$Then $\frakI_r(n)_X$ is concentrated in $[0,n]$. Moreover, the triple $(\frakI_r(n)_X, t, g)$ is unique up to a unique isomorphism and $g$ is determined by $(\frakI_r(n)_X, t)$.
\end{defn}

The duality theorem induced by the complex $\frakI_r(n)_X$ is as follows
\begin{thm}[Sato, \mbox{\cite[Theorem 1.2.2]{sato}}]
Assume that $D$ is an number ring and $X$ is proper over $S=\Spec
D$ satisfying Condition (\ref{condition*}). $\dim_SX=d$. Then there is a perfect
pairing of finite groups:
$$H^q_c(X,\frakI_r(n)_X)\times H_{\text{\'{e}}t}^{2d+1-q}(X,\frakI_r(d-n)_X) \longrightarrow H^{2d+1}_c(X,\frakI_r(d)_X)\stackrel{\cong}\longrightarrow \Z/p^r.$$
\end{thm}

The following proposition and theorem compare $\frakI_r(n)_X$ with Bloch's cycle
complex, which partially proves a conjecture made by Sato \cite[Conjecture 1.4.1]{sato}. If $X$ is smooth over $D$, the conclusion (without truncation) was proved by Geisser in \cite[Theorem 1.3]{geisser3}.
\begin{prop}\label{prop:cycletosato}
Let $X$ satisfy Condition (\ref{condition*}), then

(1) There is an isomorphism in $D^b(X_{\text{\text{\'{e}}t}},
\Z/p^r)$
$$\frakI_r(n)_X  \stackrel{\cong}\longrightarrow \tau_{\le n}(\Z^c_X/p^r(d-n)[-2d]).$$

(2)There is an isomorphism in $D^b(X_{\text{Zar}}, \Z/p^r)$
$$  \tau_{\le n}R\epsilon_*\frakI_r(n)_X \stackrel{\cong}\longrightarrow \tau_{\le n}(\Z^c_X/p^r{(d-n)}^{\text{{Zar}}}[-2d]).$$
\end{prop}
Recall that we use $\Z^c_X/p^r(d-n)^{\text{Zar}}$ to denote Bloch's complex of
Zariski sheaves, and $\Z^c_X/p^r(d-n)$ to denote that of \'{e}tale sheaves. Before proving the proposition, we
prove the following lemmas:

\begin{lemma}\label{lemma:HCGtosato}
Let $Z$ be a normal crossing variety over a perfect field of
characteristic $p$ and $\dim Z=d-1$,  then $H_{d+n-1}(\Z^c/p^r(n))$ admites Gersten resolution, i.e., the following sequence  of \'{e}tale sheaves is exact:
$$\xymatrix{0 \ar[r] & H_{d+n-1}(\Z^c_Z/p^r(n)) \ar[r]  & \displaystyle\bigoplus_{Z_{(d-1)}}R^{1-d-n}i_{z*}\Z^c_z/p^r(n) \ar[r]  &...   }$$                 $$\xymatrix{\ar[r]&\displaystyle\bigoplus_{Z_{(n+1)}}R^{-2n-1}i_{z*}\Z^c_z/p^r(n) \ar[r]&    \displaystyle\bigoplus_{Z_{(n)}} R^{-2n}i_{z*}\Z^c_z/p^r(n)\ar[r]& 0.}$$
\end{lemma}
\begin{proof}
   The above complex (without the first two terms) is $(C_Z^{HC}(n)/p^r)^{\acute{e}t}$, the Gersten complex of $H_{d+n-1}(\Z^c_Z/p^r(n))$. By Lemma \ref{Lemma:HCGtoMK} and Theorem \ref{thm:MKtoLRW}, we know that $(C_Z^{\text{HC}}(n)/p^r)^{\text{\'{e}}t}$
is isomorphic to $\widetilde{\nu}_{r,Z}(n)$. For normal crossing variety
$Z$, the latter is quasi-isomorphic to $\nu^{d-n-1}_{Z,r}[d+n-1]$ by
Sato, \cite[Corollary 2.5.2]{sato}. Hence,
$(C_Z^{\text{HC}}(n)/p^r)^{\text{\'{e}}t}$ is exact except at degree
$1-d-n$. Consider the niveau spectral sequence of higher Chow groups
in $\Z/p^r$-coefficients in (\ref{definitionspectralsequence}), we see that the cohomology of
$(C_Z^{\text{HC}}(n)/p^r)^{\text{\'{e}}t}$ at degree $1-d-n$ is
exactly $H_{d+n-1}(\Z^c_Z/p^r(n))$.

\end{proof}

\begin{lemma}\label{lemma2:cycletosato}Let $F$ be a discrete valuation field  of characteristic 0 with perfect residue field $k$ of characteristic  $p$. Then the following square is commutative: $$\xymatrix{CH^m(F,m,\Z/p^r) \ar[r]^{\delta}\ar[d]^-{\kappa} & CH^{m-1}(k, m-1,\Z/p^r)\ar[d]^-{\psi}\\
         H^m(F, \mu^{\otimes m}_{p^r}) \ar[r]^-{\partial} & \nu^{m-1}_{k,r}(k).}$$
         Here the map $\kappa$ is the map defined in \cite{bloch}, and $\delta$ is the defined in Section 2.  The map $\psi$ is defined as the composition:
$$\xymatrix{CH^{m-1}(k,m-1, \Z/p^r)\ar[r] & K_{m-1}^M(k)/p^r \ar[r]^-{d\log}_{\cong} & \nu_{k,r}^{m-1}(k),}$$
and $\partial$ is defined as the composition (\cite{sato}, (3.2.3))
$$\xymatrix{H^m(F,\mu^{\otimes m}_{p^r}) &
\ar[l]^-{\frakg}_{\cong} \ar[r]^-{\partial^{tame}} K_m^M(F)/p^r &
K_{m-1}^M(k)/p^r \ar[r]^-{d\log}& \nu^{m-1}_{k,r}(k).}$$
Here
$\partial^{tame}$ is the tame symbol of Milnor K-theory, $d\log$ is
the Bloch--Kato isomorphism, and $\frakg$ is Galois symbol.
\end{lemma}
\begin{proof}
From the definition of the maps, the square in the lemma is actually formed from
the following diagram
 $$\xymatrix{& CH^m(F,m,\Z/p^r) \ar[dl]^{\kappa}\ar[r]^{\delta}\ar[d]^{\chi} & CH^{m-1}(k, m-1,\Z/p^r)\ar[d]\ar[rd]^{\psi} & \\
         H^m(F, \mu^{\otimes m}_{p^r}) & K_m^M(F)/p^r\ar[l]^{\cong}_{\frakg} \ar[r]^{\partial^{tame}} & K^M_{m-1}(k)/p^r \ar[r]^{d\log}_{\cong}& \nu^{m-1}_{k,r}(k).}$$
The map $\chi$ is the Nesterenko--Suslin isomorphism. So we have to show the trapezoid is commutative. The triangle on the right commutes by definition. The square in the middle commutes by Geisser--Levine, \cite{geisserlevine}, Lemma 3.2. So we only have to analyze the triangle on the left. Since all the three maps $\kappa, \chi, \frakg$ respect the products of higher Chow groups, Milnor K-theory and Galois cohomology, so it suffices to show that $\kappa=\frakg\circ \chi$ for $m=1$. But in this case, all the three groups are identified with $F^{\times}/F^{\times p^r}$, and the three maps are just identity. Hence we proved the lemma.
\end{proof}
\begin{proof}[Proof of Proposition \ref{prop:cycletosato}](1) $\Rightarrow$ (2): From Theorem \ref{thm:BLconjecturetorsion}, we see that  $$\tau_{\le n}(\Z^c_X/p^r(d-n)^{Zar}[-2d])\cong \tau_{\le n}R\epsilon_*(\Z^c_X/p^r(d-n)[-2d])\cong \tau_{\le n}R\epsilon_*(\tau_{\le n}\Z^c_X/p^r(d-n)[-2d]).$$ Hence, apply $\tau_{\le n}R\epsilon_*$ on the isomorphism in (1), we obtain the one in (2).

(1) The idea of the proof is to show that $\Z^c_X(d-n)/p^r[-2d]$
satisfies a triangle similar to the one for
$\frakI_r(n)_X$, and then show that there is a map of these two triangles which induces isomorphisms on the cohomologies  up to degree $n$.

 Let $$M_X=\Z_X^c/p^r(d-n)^{\text{Zar}}[-2d], M_U=\Z_U^{c}/p^r(d-n)^{\text{Zar}}[-2d], M_Y=\Z_Y^{c}/p^r(d-n)^{\text{Zar}}[1-2d].$$ Then by Theorem \ref{thm:localization},
there is a distinguish triangle of complexes of Zariski sheaves:
$$ M_X \to j_*M_U \stackrel{\delta}\longrightarrow \iota_*M_Y .$$
Since
$\epsilon^*$ is exact and
$\epsilon^*\Z^c_X(n)^{Zar}=\Z^c_X(n)$, we get a
distinguish triangle of complexes of etale sheaves
\begin{equation} \Z^c_X/p^r(d-n)[-2d] \to \epsilon^*j_*M_U \stackrel{\delta}\longrightarrow \epsilon^*\iota_*M_Y .\label{distinguishtriangle1}\end{equation}

Consider the following diagram:
$$\xymatrix{\iota_*\widetilde{\nu}_{Y,r}(d-n)[1-2d]  & (\iota_*\tau_{\le n}\widetilde{\nu}_{Y,r}(d-n))[1-2d]\ar[l]^{\cong}  \ar[r]^{\cong}& H^n(\iota_*\widetilde{\nu}_{Y,r}(d-n)[1-2d])[-n]\\
\epsilon^*\iota_*M_Y \ar[u]^{\cong} & \tau_{\le n}\epsilon^*\iota_*M_Y \ar[r]\ar[l] \ar[u]^{\cong}& H^n(\epsilon^*\iota_*M_Y)[-n]\ar[u]^{\cong}\\
 \epsilon^*j_*M_U\ar[u]^{\delta} \ar[d] & \ar @{} [dr]|{\sharp} \tau_{\le n} \epsilon^*j_*M_U \ar[u] \ar[l]\ar[r] \ar[d]^{\cong} & H^n(\epsilon^*j_*M_U)\ar[u]^{\delta}\ar[d]^{\cong}[-n]\\
 Rj_*\epsilon^*M_U \ar[d] & \tau_{\le n}Rj_*\epsilon^*M_U \ar[d]^{\cong}\ar[r]\ar[l] & (R^nj_*\epsilon^*M_U)[-n] \ar[d]_{\kappa}^{\cong}\\
 Rj_*\mu_{p^r}^{\otimes d} & \tau_{\le n}Rj_*\mu_{p^r}^{\otimes d}  \ar[r] \ar[l]& (R^nj_*\mu_{p^r}^{\otimes d})[-n].  }$$
 The maps on the left column are defined as follows: the first one is the map defined in Theorem \ref{thm:cycletoLRW}, which is a quasi-isomorphism (note that $\epsilon^*\iota_*=\iota_*\epsilon^*$). The second map is the boundary map of localization sequence of higher Chow groups as above. The third map $\epsilon^*j_* \to Rj_*\epsilon^*$ is induced by universal property of the $\text{\'{e}}$tale sheafification $\epsilon^*$. The fourth map is defined in \cite[\S12]{levine2}. The vertical maps in the second and third column are  induced by taking truncation and cohomology of maps in the first column, respectively. Hence all the squares commute.

  The maps in the first row are isomorphisms, since it is the Gersten resolution of $\nu_{Y,r}^{n-1}[-n]$ on normal crossing varieties over perfect fields of characteristic $p$ \cite[Corollary 2.2.5(1)]{sato2}. In particular, there is an isomorphism of sheaves $$\iota_*\nu_{Y,r}^{n-1}\stackrel{\cong}\longrightarrow H^n(\iota_*\widetilde{\nu}_{Y,r}(d-n)[1-2d]),$$ which, composing with the isomorphism $$H^n(\epsilon^*\iota_*M_Y)\stackrel{\cong}\longrightarrow H^n(\iota_*\widetilde{\nu}_{Y,r}(d-n)[1-2d]),$$ induces an isomorphism $$\psi:H^n(\epsilon^*\iota_*M_Y) \stackrel{\cong}\longrightarrow \iota_*\nu^{n-1}_{Y,r}.$$  On the other hand, the Beilinson--Lichtenbaum Conjecture (Theorem \ref{thm:BLconjecturetorsion}) claims that
 $$\tau_{\le n}R\epsilon_*(\epsilon^*M_U)\cong \tau_{\le n}M_U.$$ 
 Therefore,  $$ \tau_{\le n}\epsilon^*j_*M_U\cong \tau_{\le n}Rj_*\epsilon^*M_U,$$
  i.e., in the square $\sharp$, the vertical maps are isomorphisms. At last, the vertical maps in the lower right
  square are isomorphisms, by Theorem \ref{thm:kummerisomorphism}.

In conclusion, we have the following diagram
\begin{equation}\xymatrix{\Z^c_X/p^r(d-n)[-2d] \ar[r] & Rj_*\epsilon^*M_U \ar[rr] &  &  \epsilon^*\iota_*M_Y \\
 & \tau_{\le n}Rj_*\epsilon^*M_U \ar[r]\ar[u] \ar[d]^-{\cong}& (R^nj_*\epsilon^*M_U) [-n] \ar[r]^{\delta}\ar[d]_{\kappa}^{\cong} \ar @{} [dr]|{\diamondsuit}& H^n(\epsilon^*\iota_*M_Y)[-n] \ar[d]_{\psi}^{\cong}\ar[u]^{\cong}\\
           \frakI_r(n)_X \ar[r] &      \tau_{\le n}Rj_*\mu_{p^r}^{\otimes n} \ar[r] &       (R^nj_*\mu_{p^r}^{\otimes n})[-n]\ar[r]^{\partial} & \iota_* \nu_{Y,r}^{n-1}[-n].}
\label{section3keysquare2}     \end{equation}
Here $\partial$ is the map $\sigma^n_{X,r}$ defined above. The top row and bottom row (with $(R^nj_*\mu_{p^r}^{\otimes n})[-n]$ eliminated) are distinguished triangles. The upper square and lower left square commute as above. If we can show that the square $\diamondsuit$ commutes, then it implies that there is a map $$\frakI_r(n)_X\to \Z^c_X/p^r(d-n)[-2d]$$
 which induces isomorphisms on cohomologies of  $\Z^c_X/p^r(d-n)[-2d]$ and $\frakI_r(n)_X$ up to degree $n$, hence finishes the proof.

Now let us prove that the square $\diamondsuit$ commutes. For $x\in U$, the
right two terms vanish, hence the square $\diamondsuit$ commutes.
 If $x\in Y$, we can take the henselization $S^{h}_{s}$ of $S$ at $s$, and prove the commutativity for all the irreducible components of $X_s=X\times_SS^{h}_s$ over $S^h_s$ (or replace $S^{h}_s$ by its integral closure in an irreducible component of $X_s$, depending on the Condition (\ref{condition*})).

   From now on, we assume that $A$ is a henselian discrete valuation ring with field of fraction $K$ of characteristic 0 and perfect residue field $k$ of characteristic $p$, and $X$ is flat, regular semistable family over $A$. Then $V=X_K$ is smooth over $K$ and $Z=X_k\stackrel{\iota}\longrightarrow X$ is a normal crossing variety over $k$. Since we are only left to show the commutativity of square $\diamondsuit$ in the diagram (\ref{section3keysquare2}) for points in $Z$, we can apply $\iota^*$ on the square, and obtain the following diagram:
\begin{equation}
\xymatrix{\iota^*R^nj_*\epsilon^*M_V \ar[r]^{\delta} \ar[d]^{\kappa}& H^n(\epsilon^*M_Z) \ar[d]^{\psi}\\
           \iota^*R^nj_*\mu^{\otimes n}_{p^r} \ar[r]^{\partial} & \nu^{n-1}_{Z,r}.}\label{section3keysquare3}
           \end{equation}
To show that this square commutes, we will use Gersten complex to
replace all the maps by maps of functions fields of the schemes, and
show compatibility of those maps.

First, we have the following isomorphism of complexes:
$$\xymatrix{
 \iota^*R^nj_*(\epsilon^*M_V) \ar[r] \ar[d]_{\cong}^{\kappa}&
\displaystyle\bigoplus_{V^{(0)}}\iota^*R^n(i_{x})_*(\epsilon^*M_x) \ar[r]\ar[d]_{\cong}
& \displaystyle\bigoplus_{V^{(1)}}\iota^*R^{n-1}(i_{x})_*(\epsilon^*M_x)\ar[d]_{\cong}
\ar[r] & \\
 \iota^*R^nj_*\mu^{\otimes n}_{p^r} \ar[r] & \displaystyle\bigoplus_{V^{(0)}}\iota^*R^n(i_{x})_*\mu^{\otimes n}_{p^r} \ar[r] & \displaystyle\bigoplus_{V^{(1)}}\iota^*R^{n-1}(i_{x})_*\mu^{\otimes (n-1)}_{p^r} \ar[r]  & .}$$
Here $i_x:x\to V \to X$ is the embedding into $X$, and $M_x=\Z^c_x/p^r(d-n)^{Zar}[-2i]$ if $x\in V_{(i)}$. Secondly, by \cite{sato2},
$\nu^{n-1}_{Z,r}$ on normal crossing variety admits Gersten
resolution. By Lemma \ref{lemma:HCGtosato}, higher Chow groups on
normal crossing varieties admits Gersten resolution as well. More
precisely, we have the following commutative diagram (\ref{thm:cycletoLRW}):
$$\xymatrix{
0 \ar[r] & H^n(\epsilon^*M_Z) \ar[r] \ar[d]^{\psi}_{\cong}&
\displaystyle\bigoplus_{Z^{(0)}}R^ni_{z*}(\epsilon^*M_z) \ar[r] \ar[d]_{\cong}&
\displaystyle\bigoplus_{Z^{(1)}}
R^{n-1}i_{z*}(\epsilon^*M_z)\ar[r]\ar[d]_{\cong}&\\
0 \ar[r] & \nu^{n-1}_{Z,r} \ar[r] & \displaystyle\bigoplus_{Z^{(0)}}i_{z*}\nu_{r,z}^{n-1} \ar[r] & \displaystyle\bigoplus_{Z^{(1)}} i_{z*}\nu_{r,z}^{n-2} \ar[r]&
}$$
such that the rows are exact. Here $i_z:z\to Z$ are points in $Z$, $M_z=\Z^c_z/p^r(d-n)^{Zar}[-2i-1]$ on $\Spec k(z)$ with $z\in Z^{(d-1-i)}=Z_{(i)}$.  In particular, the map $$\nu^{n-1}_{Z,r}\to \bigoplus_{Z^{(0)}}i_{z*}\nu^{n-1}_{r,z}$$
 is injective. Thirdly, by functoriality of cycle complex, the following diagram is commutative,
$$\xymatrix{
 & \iota^*R^nj_*(\epsilon^*M_V) \ar[r] \ar[d]^{\delta}&
\displaystyle\bigoplus_{V^{(0)}}\iota^*R^ni_{x*}(\epsilon^*M_x) \ar[r]
\ar[d]&
\displaystyle\bigoplus_{V^{(1)}}\iota^*R^{n-1}i_{x*}(\epsilon^*M_x)
\ar[d]
\\0 \ar[r] & H^n(\epsilon^*M_Z) \ar[r] & \displaystyle\bigoplus_{Z^{(0)}}R^ni_{z*}(\epsilon^*M_z) \ar[r] & \displaystyle\bigoplus_{Z^{(1)}} R^{n-1}i_{z*}(\epsilon^*M_z).}$$
Finally, the map $\partial$ is induced by the  diagram in
(\ref{definition:sigma}).

To summarize the above analysis, let $z\in Z^{(0)}$. Replacing $X$ by the component containing $z$, we can assume that $X$ is irreducible, and $i_{\eta}:\eta\to X$ is the embedding of the generic point. To prove the commutativity of square (\ref{section3keysquare3}), it suffices to prove that the following square commutes (here $i_z:z\to Z$):
$$\xymatrix{(\iota^*R^ni_{\eta*}(\epsilon^*M_x))_{\bar{z}} \ar[d]\ar[r] & (R^ni_{z*}\epsilon^*M_z)_{\bar{z}} \ar[d] \\
          (\iota^*R^ni_{\eta*}\mu^{\otimes n}_{p^r})_{\bar{z}}    \ar[r]      & (i_{z*}\nu^{n-1}_{r,z})_{\bar{z}} .}$$
Since
$$(\iota^*R^ni_{\eta*}\epsilon^*M_x)_{\bar{z}}=H^n_{\acute{e}t}(k(O^{sh}_{X,\bar{z}}),\Z^c/p^r(d-n)[-2d])
\stackrel{\natural}\cong H_{2d-n}^{\BM}(k(O^{sh}_{X,\bar{z}}), \Z/p^r(d-n)),$$ $$(\iota^*R^ni_{\eta*}\mu^{\otimes
n}_{p^r})_{\bar{z}}=H^n_{\acute{e}t}(k(O^{sh}_{X,\bar{z}}), \mu^{\otimes
n}_{p^r}),$$
$$(R^ni_{z*}\epsilon^*M_z)_{\bar{z}}=H^{n}_{\acute{e}t}(k(\bar{z}),
\Z^c/p^r(d-n)[1-2d])\stackrel{\natural}\cong H_{2d-n-1}^{\BM}(k(\bar{z}), \Z^c/p^r(d-n)),$$
$$(i_{z*}\nu^{n-1}_{r,z})_{\bar{z}} =\nu^{n-1}_{r,z}(k(\bar{z})).$$
Here $k(O^{sh}_{X,\bar{z}})$ is the field of fraction of
$O^{sh}_{X,\bar{z}}$, and the isomorphisms $\natural$ follow from the Beilinson--Lichtenbaum Conjecture (Theorem \ref{thm:BLconjecturetorsion}). Note that $$H_{2d-n}^{\BM}(k(O^{sh}_{X,\bar{z}}), \Z/p^r(d-n))=CH^n(k(O^{sh}_{X, \bar{z}}), n,
Z/p^r),$$  $$H_{2d-n-1}^{\BM}(k(\bar{z}), \Z^c/p^r(d-n))=CH^{n-1}(k(\bar{z}),n-1,\Z/p^r).$$ From Lemma \ref{lemma2:cycletosato}, we obtain the commutativity.
\end{proof}

\begin{thm} For $X$ satisfying Sato's condition and $n\in [0,d]$,\label{thm:cycletosato2} if the conjecture $\mathcal{B}(n)$ with $\Z/p^r$-coefficients holds for all the points in $X$, then
$$\frakI_r(n)_X  \stackrel{\cong}\longrightarrow \Z^c_X/p^r(d-n)[-2d],$$
$$  \tau_{\le n}R\epsilon_*\frakI_r(n)_X \stackrel{\cong}\longrightarrow \Z^c_X/p^r{(d-n)}^{\text{{Zar}}}[-2d].$$
\end{thm}
\begin{proof}
From Lemma \ref{lemma:gerstenconjecture}, we see that the cycle complexes in the right hand side of the isomorphisms are acyclic at degree $>n$. Therefore, in the isomorphisms of Theorem \ref{prop:cycletosato}, the truncation in front of cycle complexes are not necessary.
\end{proof}

\noindent\textit{Changlong Zhong}

\noindent University of Southern California

\noindent Department of Mathematics,

\noindent 3620 S Vermont Ave, KAP 108

\noindent Los Angeles, CA 90007

\noindent USA

\noindent Email: \textit{zhongusc@gmail.com}
\end{document}